\documentclass[a4paper]{amsart}

\usepackage{psfrag}

\usepackage[T1]{fontenc}
\usepackage{amssymb}
\usepackage{amsmath}
\usepackage[dvipsnames]{xcolor}
\usepackage{braket,comment,soul}
\usepackage{hyperref}
\usepackage{pdfpages,faktor}
\usepackage{graphicx,mathabx,bbm}
\usepackage{makecell}
\usepackage{caption}
\usepackage{subcaption}
\usepackage{mathrsfs} 
\usepackage{tikz-cd} 
\usepackage{bm}
\usepackage{enumitem}
\usepackage{mathtools}
\usepackage{amsfonts}
\usepackage{epstopdf}
\usepackage{algorithmic}
\usepackage{multirow}
\ifpdf
  \DeclareGraphicsExtensions{.eps,.pdf,.png,.jpg}
\else
  \DeclareGraphicsExtensions{.eps}
\fi

\newcommand{\R}{\mathbb{R}}

\newcommand{\C}{\mathbb{C}}

\newcommand{\Z}{\mathbb{Z}}

\newcommand{\D}{\mathbb{D}}
\newcommand{\cF}{\mathcal{F}}
\newcommand{\cA}{\mathcal{A}}
\newcommand{\cN}{\mathcal{N}}
\newcommand{\cR}{\mathcal{R}}
\newcommand{\cS}{\mathcal{S}}
\newcommand{\cK}{\mathcal{K}}

\newcommand{\cD}{\mathcal{D}}
\newcommand{\cH}{\mathcal{H}}
\newcommand{\cT}{\mathcal{T}}
\newcommand{\cU}{\mathcal{U}}
\newcommand{\cM}{\mathcal{M}}

\DeclareMathOperator{\PSL}{PSL}

\DeclareMathOperator{\Int}{Int}

\newtheorem{lemma}{Lemma}[section]
\newtheorem{theorem}[lemma]{Theorem}

\theoremstyle{definition}
\newtheorem{definition}[lemma]{Definition}

\usepackage{amsopn}

\begin{document}

\title[Where rational dynamics meets Kleinian groups]{Algebraic correspondences and Schwarz reflections: Where rational dynamics meets Kleinian groups}
    
\begin{author}[L.~Lomonaco]{Luna Lomonaco}
\address{Instituto de Matem{\'a}tica Pura e Aplicada, Estrada Dona Castorina 110, Jardim Bot{\^a}nico, Rio de Janeiro, RJ, CEP 22460-320, Brazil}
\email{luna@impa.br}
\thanks{L.L  was partially supported by the Serrapilheira Institute (grant number Serra-1811-26166), the FAPERJ - Fundação Carlos Chagas Filho de Amparo à Pesquisa do Estado do Rio de Janeiro (grant number JCNE - E26/201.279/2022 and JCM - E-26/210.016/2024), the ICTP through the Associates Programme and from the Simons Foundation through grant number 284558FY19}
\end{author}

\begin{author}[S.~Mukherjee]{Sabyasachi Mukherjee}
\address{School of Mathematics, Tata Institute of Fundamental Research, 1 Homi Bhabha Road, Mumbai 400005, India}
\email{sabya@math.tifr.res.in}
\thanks{S.M. was partially supported by the Department of Atomic Energy, Government of India, under project no.12-R\&D-TFR-5.01-0500, an endowment of the Infosys Foundation, and SERB research project grant MTR/2022/000248.}
\end{author}

\maketitle

\begin{abstract} 
We present an overview of the rapidly evolving field of dynamics of algebraic correspondences, with a focus on matings between rational maps and Kleinian groups. These correspondences exhibit rich dynamics, both within the Sullivan dictionary and beyond. We highlight unifying structures in their parameter spaces, showing how moduli spaces of rational maps and Kleinian groups naturally connect. We also outline a range of applications of the techniques developed in this framework and conclude with several promising directions.
\end{abstract}

\section{Introduction.}\label{intro_sec}

Several analogies, some philosophical and others more explicit, exist between two branches of conformal dynamics: iterations of rational maps and actions of
Kleinian/Fuchsian groups on the Riemann sphere. While the study of rational dynamics started with groundbreaking works of Fatou and Julia in the 1920s, Fuchsian/Kleinian groups were first studied by Poincar{\'e}, Fuchs, and Klein decades earlier.

Under the action of a rational map $R$, the Riemann sphere $\widehat \C$ splits into two completely invariant subsets: the (open) \textit{Fatou set} $\cF_R$, where the family of iterates $R^n$ is equicontinuous, or equivalently normal, and the (compact) Julia set $J_R:=\widehat \C \setminus \cF_R$, where the dynamics is chaotic. If the Julia set $J_R$ is not the whole Riemann sphere, then it has empty interior, and if $\deg{(R)}\geq 2$, then $J_R \neq \emptyset$, is a perfect set, is the closure of repelling periodic orbits, and is the smallest closed invariant set with at least $3$ points (cf. \cite{Bea91,CG93,Mil06}).

A non-elementary Kleinian group $\Gamma$; i.e., a discrete subgroup of $\PSL(2,\C)$ that is not virtually abelian, 
gives a dynamically invariant partition of $\widehat\C$ into two sets: the \textit{ordinary set} (or \emph{domain of discontinuity}) $\Omega_{\Gamma}$, the largest open set where the elements of $\Gamma$ form a normal family (or $\Gamma$ acts properly discontinuously), and the (compact, perfect) \textit{limit set} $\Lambda_{\Gamma}=\widehat{\C}\setminus\Omega_\Gamma$ that is the set of accumulation points of $\Gamma$-orbits. The limit set is the closure of fixed points of loxodromic elements. Further, if $\Lambda_{\Gamma}\neq\widehat{\C}$, it has empty interior (cf. \cite{Bea95,Mas88,Mar16}). 

Guided by these similarities, Fatou suggested that a comprehensive understanding of the dynamics of algebraic correspondences (which contain the above objects as sub-classes) may provide a unified framework for the dynamical study of rational maps and Kleinian groups \cite{Fat29}. But there are serious obstacles in developing a general dynamical theory for correspondences \cite{BP01}. In his landmark paper \cite{Sul85}, Sullivan introduced quasiconformal methods in complex dynamics and proposed a dictionary between rational maps and Kleinian groups. Several efforts to add precise connections in the \emph{Sullivan dictionary} were made in the next decades  \cite{McM91,MS98,LM97,Pil03}.

In this expository article, we attempt to give a unified account of some recent progress towards the realization of Fatou's original vision. An important question on this theme is the following: can these two worlds be combined in a systematic way; i.e., does there exist a mechanism to produce objects that exhibit dynamical features of rational maps and Kleinian group simultaneously? Such an object would be called a \textit{mating}.
We start by presenting matings between rational maps and Kleinian groups within holomorphic correspondences (\S\ref{BP_sec}). We then discuss how antiholomorphic counterparts of such matings arise from a different perspective involving Schwarz reflections (\S\ref{schwarz_sec}). 
We formulate general combination recipes (\S\ref{gen_matings_sec}), and talk about parameter space implications of such combination theorems (\S\ref{para_unif_subsec}); namely, the co-existence of rational parameter spaces and Kleinian deformation spaces
in parameter spaces of algebraic correspondences. 
The combination framework expounded in this article is an extension of classical combination results in groups, geometry, and dynamics; such as the \emph{Klein combination theorem} \cite{Mas93}, \emph{Bers simultaneous uniformization theorem} \cite{Ber60}, \emph{Thurston's double limit theorem} \cite{Thu86b}, and \emph{matings between polynomials} \cite{Dou83,Tan92} in a hybrid setting (cf. \cite{Bul00}). We also collect various new lines in the dictionary motivated by mating constructions, and mention how these links allow for translation of results from one world to the other (\S\ref{limit_julia_homeo_sec}, \S\ref{applications_sec}). We end with a panorama of related works (\S\ref{connections_sec}) and open problems (\S\ref{new_directions_sec}).

\noindent\textbf{Notation.} $\D=\{z\in\C:\vert z\vert<1\}$. $X^\complement=\widehat{\C}\setminus X$. $C_n$= The cyclic group of order $n$. $\eta^-(z):=1/\overline{z}$.
\section{The first family of matings.}\label{BP_sec}
In order to find a mating between a rational map and a Kleinian group, we first need to find an object that can behave like a rational map and as a Kleinian group, and then show that the dynamics of some map is compatible to the dynamics of some group. The first task brings us to the world of \textit{holomorphic correspondences}.
\smallskip

\noindent\textbf{Holomorphic correspondences.}
An $(n:m)$ holomorphic correspondence on $\widehat{\C}$ is a multivalued map $\cF: z \rightarrow w$ defined by $P(z,w)=0$, where $P(z,w)$ is a complex polynomial of degree $n$ in $z$ and degree $m$ in $w$.
A degree $d$ rational map $R:=\frac{P}{Q}$ can be written as an $(n:1)$ holomorphic correspondence by taking $P(z,w)=P(z)-wQ(z)$. A Kleinian group $\Gamma$ generated by $g_j(z)=\frac{a_jz+b_j}{c_jz+d_j},\,\,j \in [1,d],$ can be written as a $(d:d)$ holomorphic correspondence by taking $P(z,w)= \prod_{j=1}^n (w(c_jz+d_j)-(a_jz+b_j))$.
Let $R$ be a degree $d$ rational map, its \textit{deleted covering correspondence} is the $(d-1:d-1)$ correspondence $\cF$ given by the relation $\left(R(w)-R(z)\right)/(w-z)=0$.
\smallskip

\noindent\textbf{Quadratic polynomials on the Riemann sphere.}
Let $P_c(z)=z^2+c,\,\,c \in \C$. For $c \in \C$ (as for any polynomial of degree $d\geq 2$ on $\widehat{\C})$, the point at $\infty$ is fixed and superattracting, with basin of attraction $\cA_c(\infty)$. The \textit{filled Julia set $K_c=K(P_c)$} is the complement of this basin of attraction: $K_c:=\widehat\C \setminus \cA_c(\infty)$, and the Julia set $J_c=J(P_c)$ is the boundary of the filled Julia set. The \textit{Mandelbrot set} $\mathcal{M}$ is the set of $c \in \C$ such that the filled Julia set $K_c$ of $P_c$ is connected. For $c\in\mathcal{M}$, there exists a conformal map $B_c:\mathcal{A}_c(\infty)\to\overline{\D}^\complement$, called the \emph{B\"ottcher map}, that conjugates $P_c$ to the model dynamics $P_0(z)=z^2$.
If $J_c$ is locally connected, the inverse $B_c^{-1}$ extends continuously to a (semi-)conjugacy between $z^2\vert_{\mathbb{S}^1}$ and $P_c\vert_{J_c}$.
Note that the Julia set $J_0$ of $P_0(z)=z^2$ is the unit circle $\mathbb{S}^1$, and that the dynamics of $P_0$ on $\mathbb{S}^1$ is the \emph{angle doubling map}. We refer the reader to \cite{Mil06,CG93,DH84,BH88} for a detailed account of dynamics of polynomials and their connectedness loci (in particular, the Mandelbrot set).
\smallskip

\begin{figure}[h!]
\captionsetup{width=0.96\linewidth}
  \centering
 \includegraphics[width=1\linewidth]{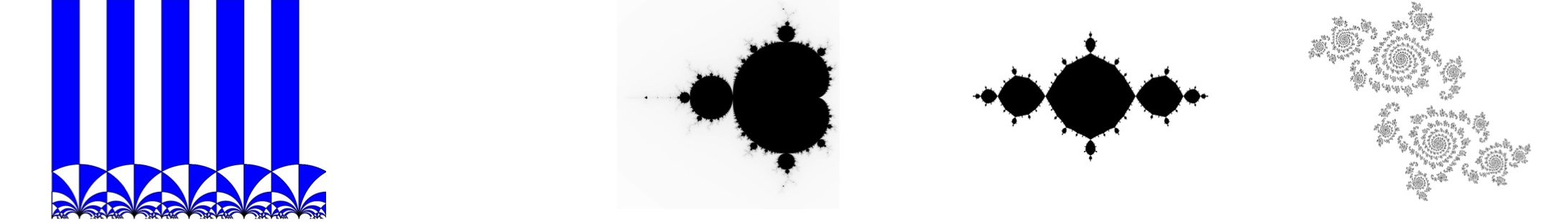}
      \caption{\small Left: Tessellation of $\PSL_2(\Z)$. Center and Right: Mandelbrot set, basilica, and a disconnected Julia set.}
  \label{blu}
\end{figure}

\noindent\textbf{The modular group.}
The modular group $\PSL(2, \Z)$ is the group of M\"obius maps
$z \mapsto \frac{az+b}{cz+d}$
with $a,b,c,d \in \Z$ and $ad-bc=1$. Many choices of
pairs of generators for the modular group exist; for instance, $\{T_1,\beta\}$, or $\{T_1,S\}$, or $\{S,\rho\}$; where
$$
T_1(z)=z+1,\quad \beta(z)=z/(1+z),\quad S(z)=-1/z,\quad \textrm{and}\quad \rho(z)= T_1 \circ S(z)=-1/(z+1).
$$
The limit set of the modular group is the extended real line $\widehat\R$. A \textit{fundamental domain} for the action of $\PSL(2, \Z)\cong C_2\ast C_3$ on the upper half-plane $\mathbb{H}$ is $D:=\{z \in \mathbb{H} \,|\, |z|>1,\,-1/2 <Re(z) <1/2\}$ (see Figure~\ref{blu}, cf. \cite[\S 9.4]{Bea95}).

\noindent\textbf{The (extended) Minkowski map.}
The \emph{extended Minkowski map} $h_+$ is a homeomorphism given by
$$h_+:[0,\infty) \to  [0,1],\quad x=[x_0; x_1, x_2, \ldots]\mapsto h_+(x)=0.\underbrace{1\ldots1}_{x_0}\underbrace{0\ldots0}_{x_1}\underbrace{1\ldots1}_{x_2}\ldots,
$$
where $[x_0; x_1, x_2, \ldots] =x_0+\cfrac{1}{x_1+\cfrac{1}{x_2+\ldots}}$ is the \textit{extended continuous fraction expansion} of $x\in \R$.
The map $h_+$ conjugates the action of $T_1(z)$ and $\beta(z)$ on $(-\infty,0]$ with the doubling 
 map, and on $[0, \infty)$ with the halving map (more precisely, there exist six commutative diagrams, see \cite[\S 7.8.2]{BF14}).

\noindent Recall that the action of $P_0(z)=z^2$ on its Julia set $J_0=\mathbb{S}^1$ is the angle doubling map. 
So, via $h_+$, the action of $P_0\vert_{\mathbb{S}^1}$ (and its inverse) fits the dynamics of the modular group on its limit set.
Hence, if $J_c$ is locally connected, then the dynamics of $P_c$ (and its inverse) on its Julia set is (semi-)conjugated by $h_+$ and the B\"ottcher map $B_c$ to the dynamics of the modular group on its limit set.
\smallskip

\noindent\textbf{Topological mating.}
The Minkowski map allows us to construct topological $2:2$ correspondences on $2$-spheres behaving like quadratic polynomials (and inverse) on a completely invariant subset of the domain, and as the modular group on its complement. The following construction is due to Bullett and Penrose \cite{BP94a}.

\noindent Let $P_c(z)=z^2+c, c\in\mathcal{M}$, be such that $J_c$ is locally connected. Let $\cK_- \vee \cK_+$ be the union
of $2$ copies of $K_c$ glued together
at the $\beta$ fixed point, and define  the  $2:2$ correspondence $H_c: \cK_- \vee \cK_+ \rightarrow \cK_- \vee \cK_+$ by sending

$$
z \in \cK_-\ \textrm{ to }\ P_c(z) \in \cK_-\ \textrm{ and } -z \in \cK_+; \textrm{ and} 
$$
$$
 z\in \cK_+\ \textrm{ to the pair of points }\ P_c^{-1}(z) \in \cK_+.
$$

\noindent Let $\mathcal{G}: z \rightarrow w$ be the $2:2$ correspondence on $\mathbb H$ 
defined by the generators $T_1,\beta$ of the modular group.
We can now glue $H_c$ with $\mathcal{G}$ via the Minkowski map, obtaining  a topological correspondence on a topological $2$-sphere.
\smallskip

\noindent\textbf{Mating quadratic polynomials with the modular group.}
In \cite{BP94a} S. Bullett and C. Penrose defined a \textit{mating between a quadratic polynomial $P_c$ with connected Julia set and the modular group $\PSL(2,\Z)$} to be a $2:2$ holomorphic correspondence $\cF: z \rightarrow w$ such that: (i) there exists a completely invariant simply connected domain $\Omega$ and a biholomorphism 
 $\phi: \Omega \rightarrow \mathbb{H}$ conjugating $\cF|_{\Omega}$ to the generators of the modular group $T_1|_{\mathbb{H}}$ and $\beta|_{\mathbb{H}}$, and (ii) $\widehat \C \setminus \Omega = \Lambda= \Lambda_- \cup \Lambda_+$, where $\Lambda_- \cap \Lambda_+$ is a singleton, such that there exist a hybrid conjugacy $\varphi_{\pm}: \Lambda_{\pm} \rightarrow K_c$ such that
 $\cF|_{\Lambda_-}$ to $P_c|_{K_c}$ and  $\cF|_{\Lambda_+}$ to $P_c^{-1}|_{K_c}$, respectively. Further, they proved that: 
\begin{theorem}
For  $a = 4$, the correspondence  $\cF_a: z \rightarrow w$ given by
\begin{equation}
\left(\frac{aw-1}{w-1}\right)^2+\left(\frac{aw-1}{w-1}\right)\left(\frac{az+1}{z+1}\right)+\left(\frac{az+1}{z+1}\right)^2=3,
\label{bp_corr_eqn}
\end{equation}
is a mating between $\PSL(2,\Z)$ and $P_{-2}(z)=z^2-2$. Further, all matings between the modular group and quadratic polynomials belong to the family $\cF_a$.
\end{theorem}

 \subsection{The family $\cF_a$.}\label{bpl_family_subsec}
\begin{figure}[htbp]
 \captionsetup{width=0.96\linewidth}
\includegraphics[width=0.96\linewidth]{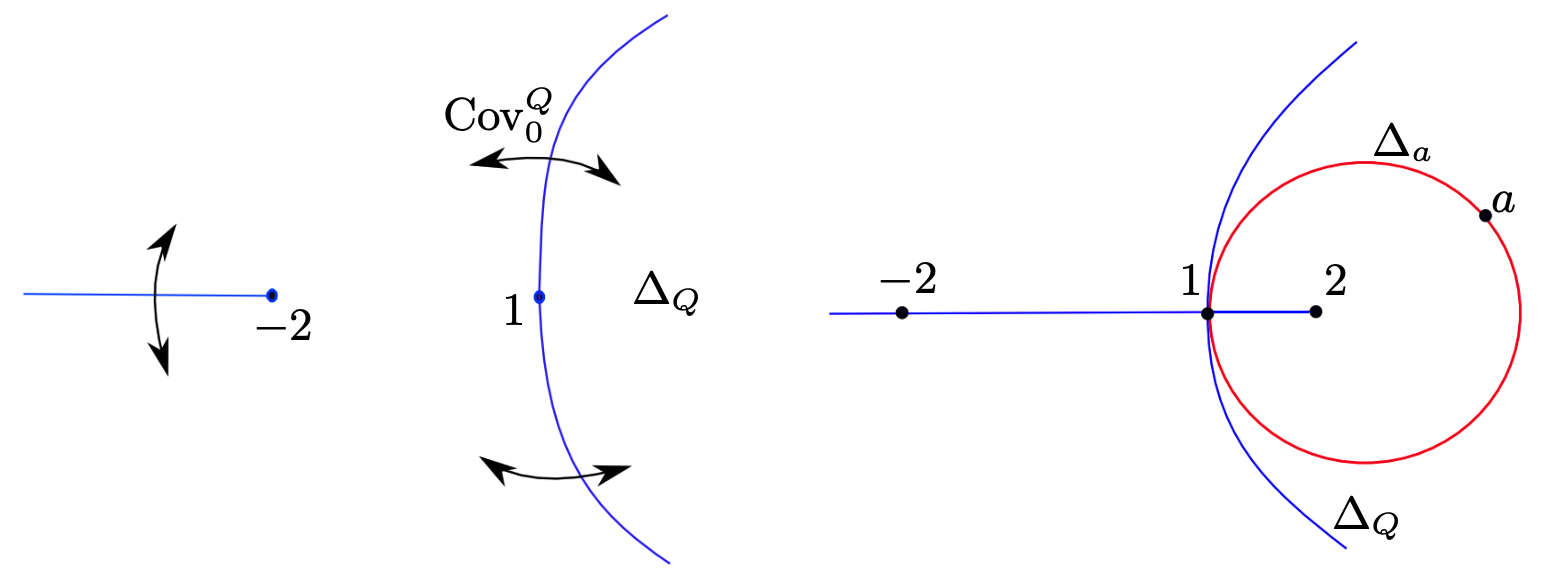}
\includegraphics[width=0.75\linewidth]{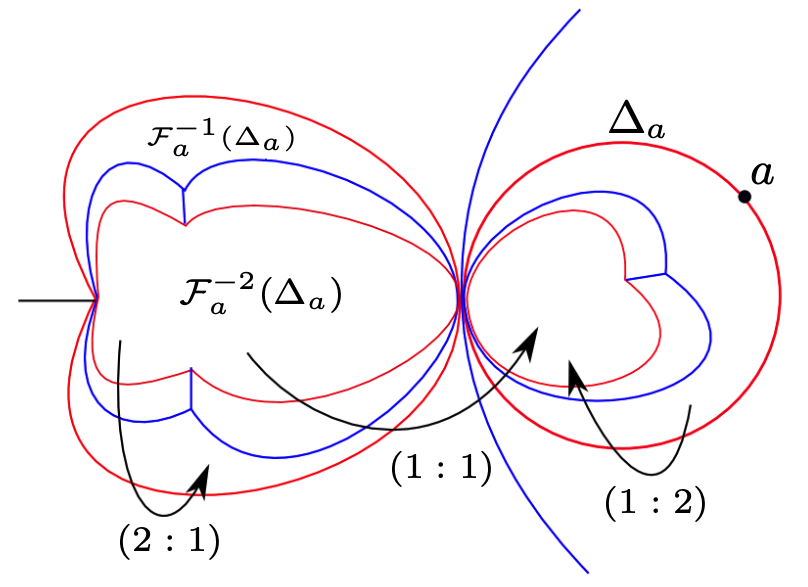}
  \caption{Top: Schematic fundamental domains $\Delta_Q$ for $Cov_0^Q$ and $\Delta_a$ for $J_a$. Bottom: Dynamics of $\cF_a$.}
  \label{ble}
\end{figure}

As a correspondence,  $\PSL(2,\Z)$ can be written as $G(z)=J \circ Cov_0^{Q_0}$, where $Q_0(z)=z^3$, $J$ is an involution, and $Cov_0^{Q_0}$ is the $2:2$ deleted covering correspondence defined in the beginning of this section.
The correspondence  $\cF_a: z \rightarrow w$ given by Equation~\eqref{bp_corr_eqn} is conjugate to
$\cF_a=J_a \circ  Cov_0^Q$, with $Q(z)=z^3-3z$, and $J_a$ being an involution with fixed point at $1$ (which is a critical point for $Q$), and $a \in \C$. 
As $Q(-2)=-2=Q(1)$, $Cov_0^Q$ sends $(-\infty,-2]$ to  curves $\ell^\pm$ emanating from $1$ and reaching $\infty$ with angle $\pm \pi/3$, so that we can take as fundamental domain $\Delta_Q$ for $Cov_0^Q$ the Jordan domain bounded by $\ell^+, \ell^-$ (see Figure \ref{ble}, left).
We consider as fundamental domain $\Delta_a$ for $J_a$ the complement of the closed round disc passing through $1, a$ having its center on the real axis (see Figure \ref{ble}, center). 
The \textit{Klein combination locus} $\mathcal{K}$ is the set of $a \in \C$ where $\cF_a$ admits fundamental domains, and $\Delta_Q^\complement \cap \Delta_a^\complement= \{1\}$, see \cite[Definition~3.1]{BL20}.
Hence, for $a\in \mathcal{K}$, $\cF_{a}: \widehat \C \setminus \Delta_a \longrightarrow \cF_a(\widehat \C \setminus \Delta_a)$ is a holomorphic $1:2$ map, while $\cF_{a}:\cF_a^{-1}(\overline \Delta_a) \longrightarrow \overline\Delta_a$ is a holomorphic $2:1$ map (see Figure \ref{ble}, right).
We define the \textit{backward limit set} $\Lambda_{a,-}$ for $\cF_a$ as $\Lambda_{a,-}=\bigcap_1^\infty(\cF_a)^{-n}(\overline \Delta_a)$, the \textit{forward limit set} $\Lambda_{a,+}$ for $\cF_a$ as $\Lambda_{a,+}= \bigcap_1^\infty (\cF_a)^n(\widehat \C \setminus \Delta_a) $ (see Figure \ref{blu}, on the top), and the \textit{limit set} $\Lambda_a:=\Lambda_{a,-} \cup \Lambda_{a,+}$.
The \textit{modular Mandelbrot set} $\cM_{\Gamma}$ is the set of parameters $a$ such that the limit set $\Lambda_a$ is connected, see Figure~\ref{blu}.
\smallskip

\noindent\textbf{Historical conjectures and challenges.}
In \cite{BP94a}, it was conjectured that $\cF_a$ contains matings between $\PSL(2,\Z)$ and $P_c$, for each $c$ in the Mandelbrot set $\cM$, and that $\cM_{\Gamma}$ is homeomorphic to $\cM$. The group $\PSL_2(\Z)$ lies on the boundary of $\cD$, which is the moduli space of discrete, faithful representations of $C_2*C_3$ in $\PSL(2,\C)$ (the modular group arises when certain loxodromic elements of these representations turn parabolic). In \cite{BH00}, matings between quadratic polynomials and representations in $\Int{\cD}$ (with Cantor limit set) were constructed via quasiconformal surgery (see \S\ref{holo_gen_mating_subsec}). One expects that appropriate limits of these matings give rise to correspondences in $\cM_\Gamma$ realizing matings between $\PSL(2,\Z)$ and $P_c$; and this was proved in \cite{BH07} for a large collection of $c \in \cM$ using pinching deformation techniques. However, the complete picture of $\cM_\Gamma$ remained out of reach at that~time.

\subsection{Matings between parabolic rational maps and the modular group.}\label{para_rat_mod_group_mating_subsec}

For every $a$ (for which the correspondence $\cF_a$ present
fundamental domains) the correspondence $\cF_a$ has a parabolic fixed point of multiplier $1$ at $z=0$ (cf. \cite[Proposition~3.5]{BL20}), which renders the strategy of \cite{BH07} inefficient. The presence of a persistent parabolic fixed point suggests the need of parabolic techniques in order to deal with the parabolic point. This was already known in \cite{BP94a} (see page 484).
In the 2010s the introduction of \textit{parabolic-like maps} by the first author \cite{Lom14,Lom15,LPS17} allowed for a systematic treatment of the family $\cF_a$, as mating between \textit{parabolic quadratic maps} and the modular group. 
More precisely, consider the family of parabolic quadratic rational maps $P_A(z)=z+1/z+A$, $A\in\C$, normalized by having a parabolic fixed point of multiplier $1$ at $\infty$ and critical points at $\pm 1$. Note that $P_A$ is conformally conjugate to $P_{-A}$; in Milnor’s notation, the set of conformal conjugacy classes is denoted by $Per_1(1)$. Let $\mathcal{A}_A(\infty)$ denote the parabolic basin of infinity of $P_A$, and note that it is completely invariant. Hence, its complement is also completely invariant, we call it the \textit{filled Julia set} $K_A:= \widehat \C \setminus \mathcal{A}_A(\infty)$ of $P_A$.  
\begin{figure}[h!] 
\captionsetup{width=0.96\linewidth}
\includegraphics[width=0.99\linewidth]{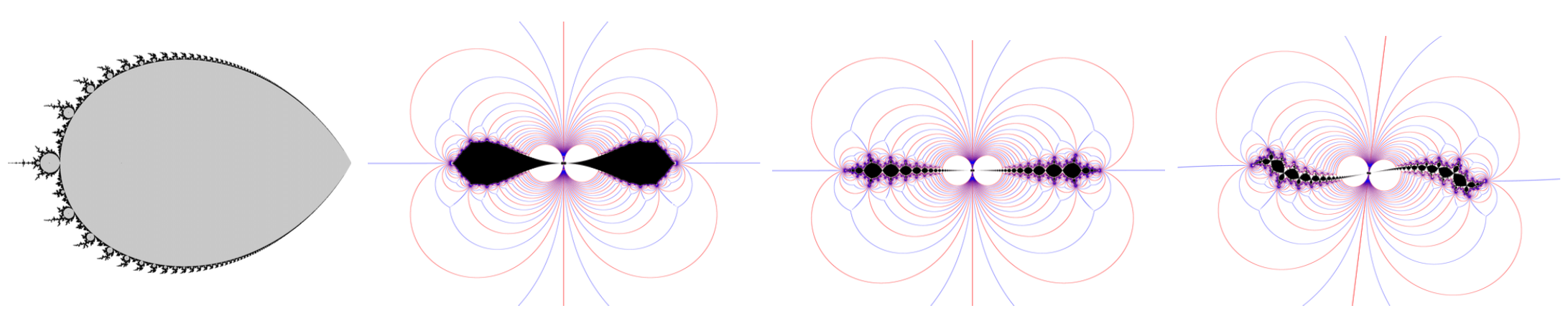}
\includegraphics[width=0.99\linewidth]{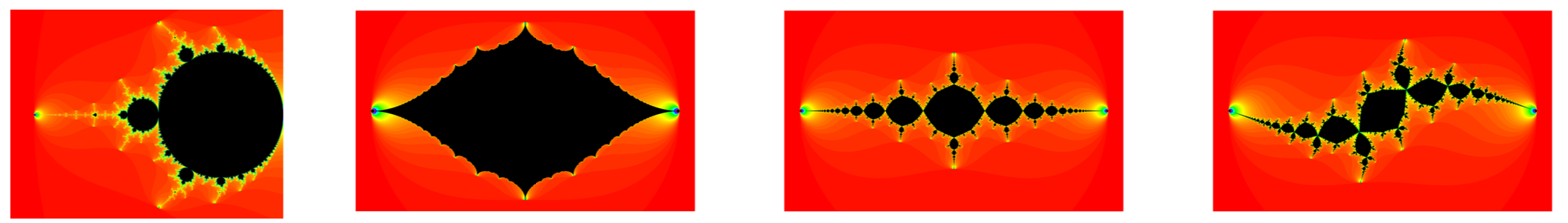} 
    \caption{Above, from the left: the modular Mandelbrot set $\cM_{\Gamma}$, and the limit sets for correspondences $\cF_a$ where $a$ is the center of the hyperbolic components of $\cM_{\Gamma}$ of periods $1, 2, 3$, respectively. Below, from the left: the parabolic Mandelbrot set $\mathcal{M}_1$, and the filled Julia sets $K_A$ for $P_A$ hybrid equivalent to $\cF_a$ above.}
  \label{fig:lambda_variations}
\end{figure}
The \textit{parabolic Mandelbrot set} is the connectedness locus for the family $Per_1(1)$, this is the set of $B=1-A^2$ such that the filled Julia set $K_A$ of $P_A$ is connected (see Figure~\ref{fig:lambda_variations}).
A parabolic-like map is an object that encodes the dynamics of $P_A$ about its filled Julia set $K_A$. More precisely, a degree $d$ parabolic-like map is a quadruple $(f,U',U,\gamma)$ where $U', U, U \cup U' \approx \D$, $U'$ is not contained into $U$, $f: U' \rightarrow U$ is a proper, holomorphic, degree $d$ map with a parabolic fixed point $P$ of multiplier $1$, and $\gamma= \gamma_- \cup \gamma_+$ are forward invariant arcs emanating from the parabolic fixed point and dividing $U$ and $U'$ in $\Delta,\,\mathcal{O}$ and $\Delta',\,\mathcal{O}'$ respectively, such that $\mathcal{O}' \subset \subset U$ and $\Delta'$ does not contain critical points (see \cite[Definition~3.1]{Lom15}). These forward invariant arcs $\gamma= \gamma_- \cup \gamma_+$ are called \textit{dividing arcs}; they separate the parabolic dynamics from the expanding one, and they meet $\partial U$ non-tangentially. Note that the dividing arcs can be taken to be (parts of) horocycles at the parabolic fixed point. The \textit{filled Julia set} for a parabolic-like map $(f,U',U,\gamma)$ is the set of points that never leave $\mathcal{O}'\cup \{P\}$ under iteration. A degree $d$ parabolic-like map is associated with an internal class, which roughly speaking is given by its filled Julia set, and its external class, which is a degree $d$ weakly expanding covering of the circle \cite{Lom15,LPS17}. By replacing the external map of a degree $2$ parabolic-like map with the Blaschke product $h(z)=\frac{z^2+1/3}{z^2/3+1}$ the first author proved that any degree $2$ parabolic-like map is hybrid equivalent to a member of the family $Per_1(1)$, and this member is unique if the filled Julia set is connected \cite{Lom15}.
In the 2010s the first author together with Shaun Bullett defined (see \cite{BL20}): 
\begin{definition}\label{modular_mating_def}
    We say that $\cF_a$ is a mating between $P_A(z)=z+1/z+A,\,A \in \C$ and $\PSL(2,\Z)$ if

1. the $2-to-1$ branch of $\cF_a$ for which $\Lambda_{a,-}$ is invariant is hybrid equivalent to $P_A$ on $K_A$, and

2. when restricted to a $(2 : 2)$ correspondence from $\Omega$ to itself, $\cF_a$ is conformally conjugate to the pair of M{\"o}bius transformations $T_1,\,\beta$ from the complex upper half plane  to itself.
\end{definition}
By turning each member of the family $\cF_a$ into a parabolic-like map, and proving that $\cF_{a|\Omega}$ is conformally conjugate to the generators of the modular group acting on $\mathbb{H}$, they prove in \cite{BL20} (see also \cite{BLS19}).

\begin{theorem}\label{bl1}
    For each $a \in \cM_{\Gamma}$, the correspondence $\cF_a$ is a mating between some rational map $P_A(z)=z+1/z+A$, $A\in\C$, and the modular group $\PSL(2,\Z)$.
\end{theorem}

\textbf{Key ideas of the proof.}
To turn each member of the family $\cF_a$ into a parabolic-like map, one starts by showing that we can always find fundamental domains $\Delta_Q$ for $Cov_0^Q$ and $\Delta_a$ for $J_a$ with boundaries smooth at $P$ and transverse to the parabolic axis, and we
show that on a neighbourhood of the preimage, call it $S$, of the parabolic fixed point $P$, the correspondence $\cF_{a}$ behaves as the square root $z \rightarrow \pm\sqrt{z}$ (see Proposition 3.8 in \cite{BL20}).
Then one 
constructs the dividing arcs $\gamma= \gamma_- \cup \gamma_+$ by taking preimages of horizontal lines under repelling Fatou coordinates of the parabolic fixed point $P$.
Hence, one shows that the $2:1$ map $\cF_{a}:\cF_a^{-1}(\overline \Delta_a) \longrightarrow \overline\Delta_a$ restricts to a $2:1$ map $\cF_{a}: V_a' \rightarrow V_a$ such that $\Lambda_{a,-} \subset V_a'$, $V_a' \subset V_a$, $\overline V_a \cap \overline V_a' =\{P\}$ the persistent parabolic fixed point, and $\overline V_a'$ makes at $P$ an angle strictly less that $\pi$ (see  \cite[Proposition~5.2, Figure~8]{BL20}). The sets $V_a'$ are of fundamental importance, and will later be called '\textbf{lunes}' in \cite{BL24}. A lune is the intersection of two discs, and it has smooth boundary except at the intersection points, where it forms an angle less that $\pi$. Technically, just the angle at $P$ needs to be less than $\pi$, so this is what has been required for $V_a'$ in \cite{BL20}.
Finally, one adds the disc-neighbourhoods of the parabolic-fixed point $P$ and its preimage $S$ to $V_a'$ where $\cF_{a}$ behaves as the square root, and performs a surgery to turn its inverse to a homeomorphism, obtaining a parabolic-like map of degree $2$ (see \cite[Theorem~B]{BL20}).
Finally, it is shown that the Riemann map $\psi:\Omega \to \mathbb{H}$
conjugates $\cF_a|_\Omega$ to an action of $C_2*C_3$ on $\mathbb{H}$; as $\Omega_a/\cF_a$ has a cusp, this is the
standard action of $PSL(2,\Z)$
(see \cite[Theorem~A]{BL20}).

\subsubsection{Dynamics of modular matings.}
A complete dynamical theory for the family $\cF_a$ is developed in \cite{BL22}. This family of correspondences, here called \textit{modular matings}, represents the simplest one-parameter case to be fully analyzed, and thus serves as a model for investigating other families of correspondences that exhibit ``discreteness'' in their action on the sphere over certain regions of parameter space. 
\smallskip

\noindent\textbf{Böttcher map.} The authors construct a \textit{Böttcher map} $\varphi_a: \Omega_a \rightarrow \mathbb{H}$; this is a conformal isomorphism conjugating the branches of $\cF_a$ on $\Omega_a$ to the action of the generators of $\PSL_2(\Z)$ on $\mathbb{H}$. The map $\varphi_a$ is used to pull back \textit{periodic geodesics}, which play in this theory the role external rays play in the classical theory of polynomials.
\smallskip

\noindent\textbf{Geodesics.} The authors first establish landing results for these geodesics in the dynamical plane. They prove that periodic geodesics always land at points of the limit set, and that every repelling fixed point of the two-to-one branch of the correspondence that fixes $\Lambda_{a,-}$ is the landing point of exactly one cycle of periodic geodesics, while repelling periodic points of higher period are landing points of one or two such cycles. This is achieved by completely new methods, based on Sturmian sequences, rather than the classical polynomial techniques.
\smallskip

\noindent\textbf{Yoccoz Inequality}. They then prove an analogue of the Yoccoz inequality, providing strong bounds on the multipliers of repelling fixed points. If such a point has combinatorial rotation number $p/q$, the logarithm of its multiplier is confined to a small disc in the right half-plane, with radius decaying on the order of $1/q^2$ (up to logarithmic factors). This sharper-than-classical inequality places severe restrictions on the dynamics near rational rotation numbers, and highlights new structural features in the correspondence setting.
\smallskip

\noindent\textbf{The Modular Mandelbrot set is contained in a lune.} Finally, the paper gives geometric constraints on the connectedness locus of parameters. They prove that the Modular Mandelbrot set $\mathcal{M}_{\Gamma}$ is contained in a lune in the parameter plane, bounded by circular arcs meeting at two real points. Moreover, they analyze the small ``limbs'' of this set corresponding to rationals with large denominators, showing that their diameters shrink to zero and that they accumulate tangentially at a specific parameter value. This yields the first concrete description of the global shape and fine structure of the parameter space for these correspondences.
\smallskip

\begin{center}
\begin{tabular}{|c | c|}
     \hline  
{\bf Quadratic correspondences ${\mathcal F}_a$} & {\bf Quadratic polynomials $P_c$}  \\ \hline

 $\varphi_a:\widehat{\mathbb C}\setminus \Lambda_a \cong {\mathbb H}$
           &  $B_c:\widehat{\mathbb C} \setminus K(P_c)\cong \widehat{\mathbb C}\setminus \overline{\mathbb D}$            \\ \hline 
 external geodesics                 & external rays                                                           \\ \hline
\vspace{0.05mm} `periodic geodesics land'         & `periodic rays land'         \\ \hline
\vtop{\hbox{\strut repelling  cycles are landing} \hbox{\strut points of periodic geodesics}}    &  \vtop{\hbox{\strut repelling  cycles are landing} \hbox{\strut points of periodic rays}}          
\\ \hline
\vspace{0.05mm} Yoccoz inequality  $\sim (\log{q})/q^2$   & Yoccoz inequality   $\sim 1/q$                                             
\\ \hline
\end{tabular} 
\end{center}
\smallskip

\noindent\textbf{Quadratic correspondences $\cF_a$ vs Quadratic polynomials $P_c$.} Although the results parallel the Douady-Hubbard theory for polynomials, the proofs are quite different, reflecting the different dynamical setting: unlike in the polynomial case, there is no superattracting fixed point at infinity to anchor the construction.  The table above summarizes the parallels of the theories of modular matings and quadratic polynomials.

\subsubsection{The modular Mandelbrot set.}
In \cite{BL24}, the first author together with Shaun Bullett proves that 
the connectedness locus $\mathcal{M}_\Gamma$ for the family of correspondences $\mathcal{F}_a$ is homeomorphic to the parabolic Mandelbrot set $\mathcal{M}_1$. More precisely, the authors construct a parameter map $\chi : \mathcal{M}_\Gamma \to \mathcal{M}_1$ which is a homeomorphism, dynamical in nature: for each $a \in \mathcal{M}_\Gamma$, the correspondence $\mathcal{F}_a$ is a mating between $PSL(2,\mathbb{Z})$ and the quadratic polynomial $P_{\chi(a)}$, and they prove the following (for a graphic representation of the Theorem, see Figure \ref{fig:lambda_variations})
\begin{theorem}
There exists a dynamical homeomorphism $\chi: \cM_{\Gamma} \rightarrow \mathcal{M}_1$ that is conformal on the interior, and extends to a pinched neighborhood of $\cM_{\Gamma}$.
\end{theorem}

The proof was first achieved by performing the surgery of \cite{BL20} in a uniform way \cite{BL20b}, and so by a two-step surgery process: first constructing families of parabolic-like maps and then applying the parabolic-like straightening theorem. While this process yields the desired result, it requires additional control at each of the two surgery steps: (i) the surgery to pass from correspondences to parabolic-like maps, and (ii) the one from parabolic-like maps to rational maps. This was subsequently streamlined by the authors in \cite{BL24}, where they condensed the two steps, adapting the surgery done for converting parabolic-like maps into rational maps to the setting directly to the family $\cF_a$, thus skipping the intermediate parabolic-like map step. This latter surgery, more direct and elegant, is inspired by the surgery developed by the first author in her thesis \cite{Lom15}.
Both approaches rely crucially on the fact that the lune $V_a'$ moves holomorphically with the parameter throughout the parameter lune, which follows from the earlier result that the Modular Mandelbrot set is contained in a parameter lune.

After constructing the map $\chi : \mathcal{M}_\Gamma \to \mathcal{M}_1$ via the surgery, standard arguments (a la Douady--Hubbard and Lyubich) establish that $\chi$ is holomorphic on $\Int{\mathcal{M}_\Gamma}$, continuous up to $\partial\mathcal{M}_\Gamma$, and extends to a pinched neighborhood of $\mathcal{M}_\Gamma$. The delicate part of the proof is surjectivity. This is particularly challenging because $\mathcal{M}_\Gamma$ is not compactly contained in the parameter space for correspondences (the Klein combination locus $\mathscr{K}$), so standard compactness arguments do not apply. To address this, the authors carefully analyze the combinatorial structure of $\mathcal{M}_\Gamma$, showing that each limb of $\mathcal{M}_\Gamma$ is homeomorphic to a corresponding limb of $\mathcal{M}_1$ (see \cite[\S 4.5]{BL24}). 

Establishing these homeomorphisms requires controlling the geometry of each limb, including how decorations and satellite components attach, and ensuring that this structure is preserved under $\chi$. This is delicate because one must track the dynamics of critical points and their forward orbits, as well as the interactions between limbs and the root points of $\mathcal{M}_\Gamma$, which are not compactly contained. By proving that each limb is correctly mapped onto a corresponding limb of $\mathcal{M}_1$, the authors establish the surjectivity of $\chi$.
Together with results on holomorphicity, continuity, and the extension to a pinched neighborhood, this surjectivity completes the proof that $\chi$ is a homeomorphism between the modular Mandelbrot set $\mathcal{M}_\Gamma$ and the parabolic Mandelbrot set $\mathcal{M}_1$.

Finally, by combining this homeomorphism with the theorem of Petersen--Roesch that $\mathcal{M}_1$ is homeomorphic to the classical Mandelbrot set $\mathcal{M}$ \cite{PR21}, the authors conclude that $\mathcal{M}_\Gamma$ itself is homeomorphic to $\mathcal{M}$. This establishes that the modular Mandelbrot set has exactly the same topological type as the classical Mandelbrot set, providing a complete description of its global parameter structure. Moreover, the result confirms that the intricate combinatorial structure of $\mathcal{M}_\Gamma$, including its limbs and parameter decorations, corresponds precisely to those in the classical setting. In particular, the result allows one to transfer insights and techniques from classical quadratic dynamics directly to the study of $\mathcal{F}_a$, and vice versa, making the modular Mandelbrot set fully accessible to the combinatorial and analytic tools developed for $\mathcal{M}$.

\subsubsection{Tessellation of the exterior of the modular Mandelbrot set.}
\begin{figure}[h!]
\captionsetup{width=0.96\linewidth}
  \centering
\includegraphics[width= 4.2cm]{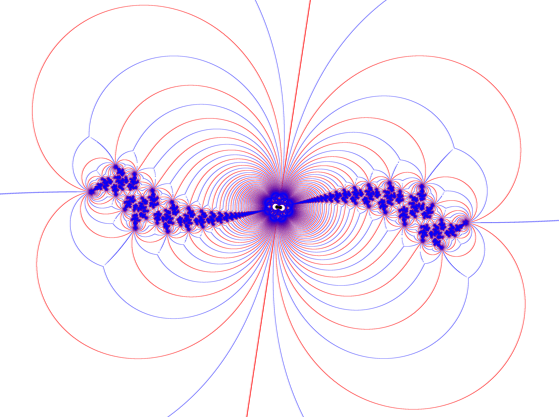}  \hspace{12mm}
\includegraphics[width= 3.6cm]{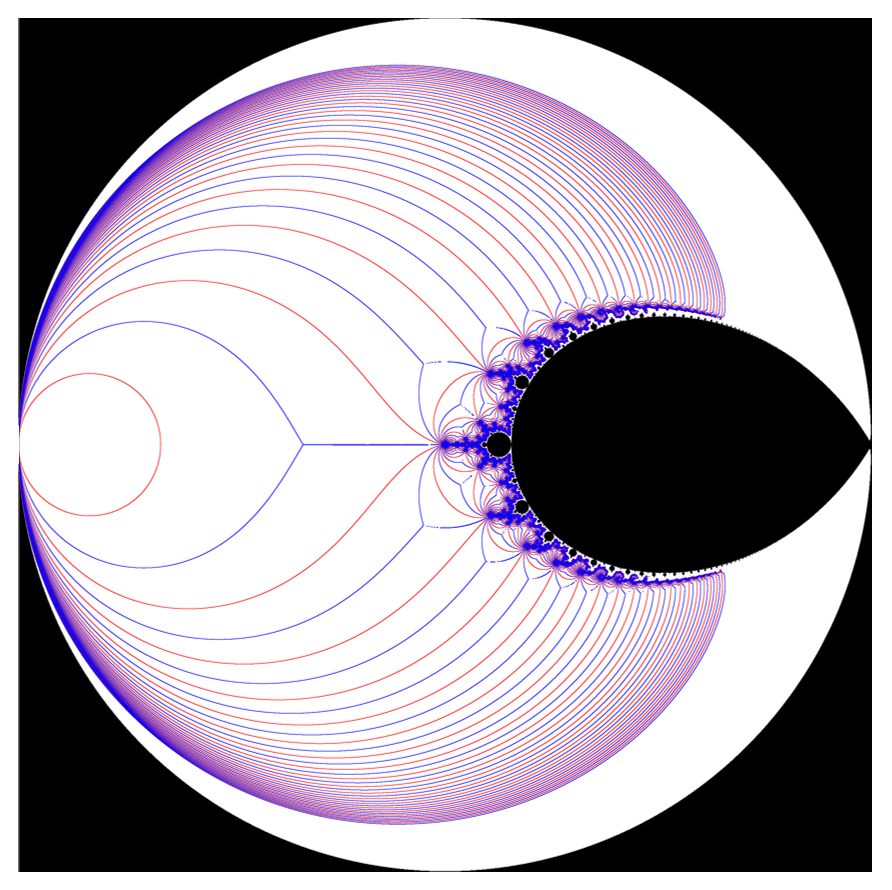}  \hspace{16mm}
\includegraphics[width= 3.2cm]{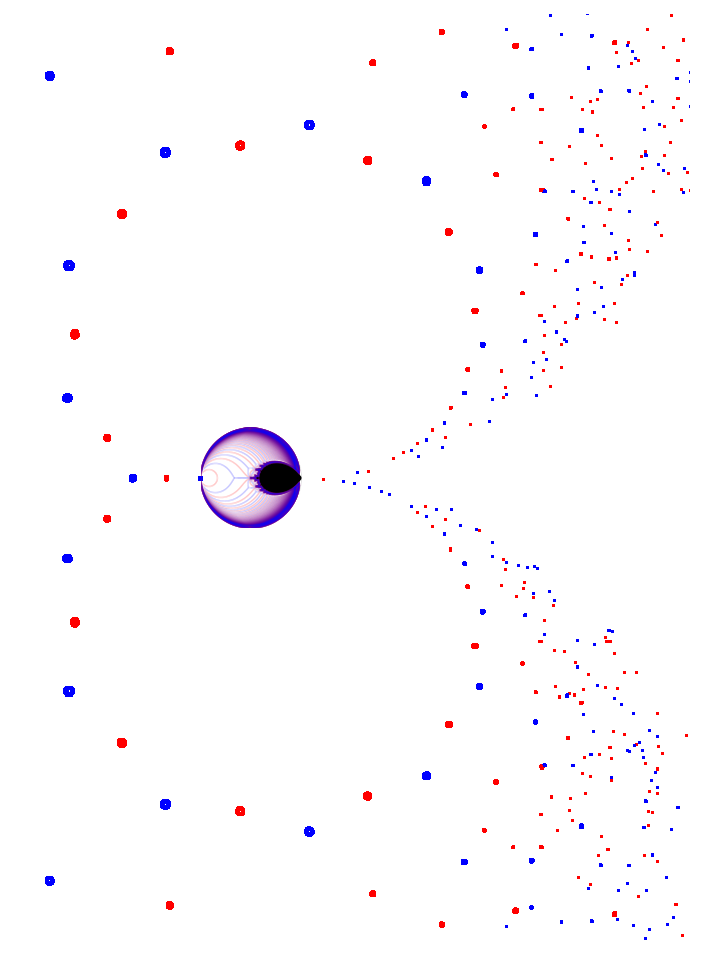}
    \caption{Top left: Dynamical tessellation via partial B\"ottcher map $\varphi_a$. Top right: Parameter tessellation on $\D(4,3)$. Bottom: Points that approximate $\partial \mathcal{K}$. Pictures are courtesy of Ivan Pedro Suarez Navarro.}
  \label{fig:tess}
\end{figure}

Recall that the \textit{Klein combination locus} $\mathcal{K}$ is the set of $a \in \C$ such that there exist fundamental domains $
\Delta_Q$ for $Cov_0^Q$ and $\Delta_a$ for $J_a$,  and $\Delta_Q^\complement \cap \Delta_a^\complement= \{1\}$. Then $\Delta_{corr}= \Delta_Q \cap \Delta_a$ is a fundamental domain for $\cF_a$. The quotient orbifolds $\Delta_{corr}/\cF_a$ and $\mathbb{H}/PSL(2,\Z)$ are conformally isomorphic \cite[Theorem A]{BL20}, so 
for all $a \in \mathcal{K}$ there exists a fundamental domain $\Delta_{mod}$ for $PSL(2,\Z)$ such that $\Delta_{corr}$ is conformally isomorphic to $\Delta_{mod}$ by say $\varphi_a$, and $\varphi_a$ sends the vertices of  $\Delta_{corr}$ to $0$, $\infty$, $i$ and $(-1+i\sqrt{3})/2$, and it is equivariant. In \cite{BLLS2}, it is shown that for every $a\in\mathcal{K}$, the critical value $v_a$ always lies in $\Delta_{corr} \cup J_a(\Delta_{corr})$ together with its inverse images. So we can extend $\varphi_a $ equivariantly until the tile containing the critical value, obtaining a partial Böttcher map $\varphi_a$ for $\mathcal{F}_a$, and define:
$$\Psi: \mathcal{K} \setminus \mathcal{M}_{\Gamma} \rightarrow \mathbb{H},\,\,\, \Psi(a)=\varphi_a(v_a)
$$

\noindent analogous to the classical construction for quadratic polynomials. Pulling back the tessellation of $\mathbb{H}$ by $\Psi$ produces a tessellation around $\mathcal{M}_\Gamma$.
For $a \in \mathbb{D}(4,3)$, $\Psi$ is a well-defined continuous injection, obtained using the “standard” $\Delta_{corr}$ (see \cite[Section 2]{BL22} and Figure \ref{fig:tess}). Changing the fundamental domains alters the tessellation, although the positions of the vertices remain fixed. In particular, encircling the puncture point $a = 1$ in the parameter plane moves the standard domains to a new location, making $\Psi$ multivalued unless we introduce a cut along $[-1,+1]$ in $\mathcal{K}$. 
The orbifold of grand orbits of $\cF_
a$ on $\Omega_a$ is a sphere with a
puncture, a $\pi$-cone point and a $2\pi/3$-cone point. The structure of $\mathcal{K}$
and its boundary in $a$-space can be at least partially determined by
classifying the values of $a$ such that the grand orbit under $\cF_
a$ and $\cF_a^{-1}$ of the critical point $c_a$ projects to one of the cone points \cite{BLLS2}.
The family of correspondences $\cF_a$ is a slice of a $2$-complex parameter family
\cite{BH00} which has as another slice the moduli space of representations of $C_2*C_3$ in $\PSL(2,\C)$. The discreteness locus
$\mathcal D$ of the latter slice (an ‘elliptic cousin of the Riley slice’ \cite{KS,EMS}) has
on its boundary a dense set of cusps, indexed by rationals $p/q$, where the
representation contains a parabolic cycle, the normality set is non-empty, and
the limit set is a circle-packing. In contrast the points $a\in \partial\mathcal K$ where $\cF_a$
has an analogous parabolic cycle, non-empty normality set and connected limit
set, are indexed by rationals of the form $1/q$, and these points of $\partial\mathcal K$
accumulate only at the root point of $\mathcal M_{\Gamma}\subset \mathcal K$. We also refer the reader to \cite{BC12} for the dynamical description of a correspondence on $\partial\cK$.

\section{A general framework of mating: internal and external maps.}\label{internal_external_sec}
Some dynamical systems on $\C$ can be partitioned into two disjoint sub-
systems: an {\it internal class} and an {\it external class}. To illustrate this, let us
consider a simple example.
For a degree $d\geq 2$ polynomial $P$ with connected filled Julia set $K(P)$, the action of $P$ on its basin of infinity $\cA_P$ is conformally conjugate to the power map $z^d$. Thus, one can regard $P$ as the union of two disjoint dynamical systems: its internal map, which is given by $P|_{K(P)}$, and its external map $z^d$ that encodes the dynamics of $P$ outside $K(P)$.
Similarly, for $R \in Per_1(1)$ with a connected $K(R)$, the dynamics of $R$ on its parabolic basin $\cA_R$ is conformally conjugate to the Blaschke product $h$ (see \S\ref{para_rat_mod_group_mating_subsec}), and hence $R$ can be thought of as the union of its internal map $R|_{K(R)}$, and its external map $h$ modeling the dynamics of $R$ outside $K(R)$.
This extends naturally to arbitrary degree; more precisely, the dynamics of a degree $d\geq 2$ rational map $R$ with a parabolic fixed point of multiplier $1$ having a connected and simply connected parabolic basin $\cA_R$ splits into its internal map $R|_{K(R)}$, and its external map given by a parabolic Blaschke product. 
A polynomial-like map is a partially defined dynamical system mimicking the action of a polynomial. More precisely, a degree $d$ polynomial-like map is a triple $(f,U',U)$ such that $U,\,U' \cong \D,\,\,\overline U' \subset U$, and $f:U' \rightarrow U$ is a degree $d$ proper holomorphic map. The action of $f$ on its \emph{non-escaping set} $K(f)$, consisting of points that never leave $U'$, is \emph{hybrid conjugate} (by a quasiconformal conjugacy $\phi$ with $\overline\partial \phi=0$ on $K(f)$) to the action of a polynomial $P$ on $K(P)$, we call $f_{|K(f)}$ the {\it internal class}. The external map of $f$ is a real-analytic, expanding circle covering capturing the dynamics of $f$ outside $K(f)$, and the {\it external class} is the real-analytic conjugacy class of the external map \cite{DH85,Lyu25}.
A parabolic-like map $(f,U',U,\gamma)$ (see \S\ref{para_rat_mod_group_mating_subsec}) also partitions into an internal class, given by the dynamics of $f\vert_{K(f)}$, and an external map, given by a real-analytic, weakly expanding circle covering encoding the dynamics of $f$ outside $K(f)$.

Matings of groups and maps can often be expressed in the language of internal and external maps. To facilitate this, we need to consider \textit{degenerate polynomial-like maps}, or \textit{degenerate p-l maps} (for brevity). Roughly speaking, these are triples 
$(f,U',U)$ where $\overline{U},\,\overline{U'}$ are (finitely) pinched disks with $\overline{U'}\subset \overline{U}$, such that $\partial U' \cap \partial U$ consists of finitely many `parabolic cycles', and $f$ is a proper holomorphic map from components of $U'$ to those of $U$ (cf. \cite{BLLM25,LLM24}). When $U\cong \D$, they are called \textit{pinched p-l maps}. A degenerate p-l map $f$ can be regarded as a mating of its internal class, given by the action of $f$ on its non-escaping/filled Julia set $K(f)$, and an external map, a piecewise real-analytic, weakly expanding circle covering that encodes the action of $f$ outside~$K(f)$. By selecting external maps from groups appropriately, and mating them with parabolic rational maps or polynomials using quasiconformal or trans-quasiconformal/David surgery, one can produce degenerate p-l maps that realize combinations of groups and polynomials.
Specifically, mating parabolic internal classes and groups with a unique conjugacy class of parabolics can be done quasiconformally \cite{BL20,LMM24,BLLM25}, as well as mating hyperbolic internal classes with groups without parabolics \cite{BH00,Lau24}; while any other mating requires trans-quasiconformal surgery \cite{LMMN25,MM25,LLM24}.
This theme will unify the combination theorems discussed next, as shown in the table.

\begin{center}
\includegraphics[width=0.8\linewidth]{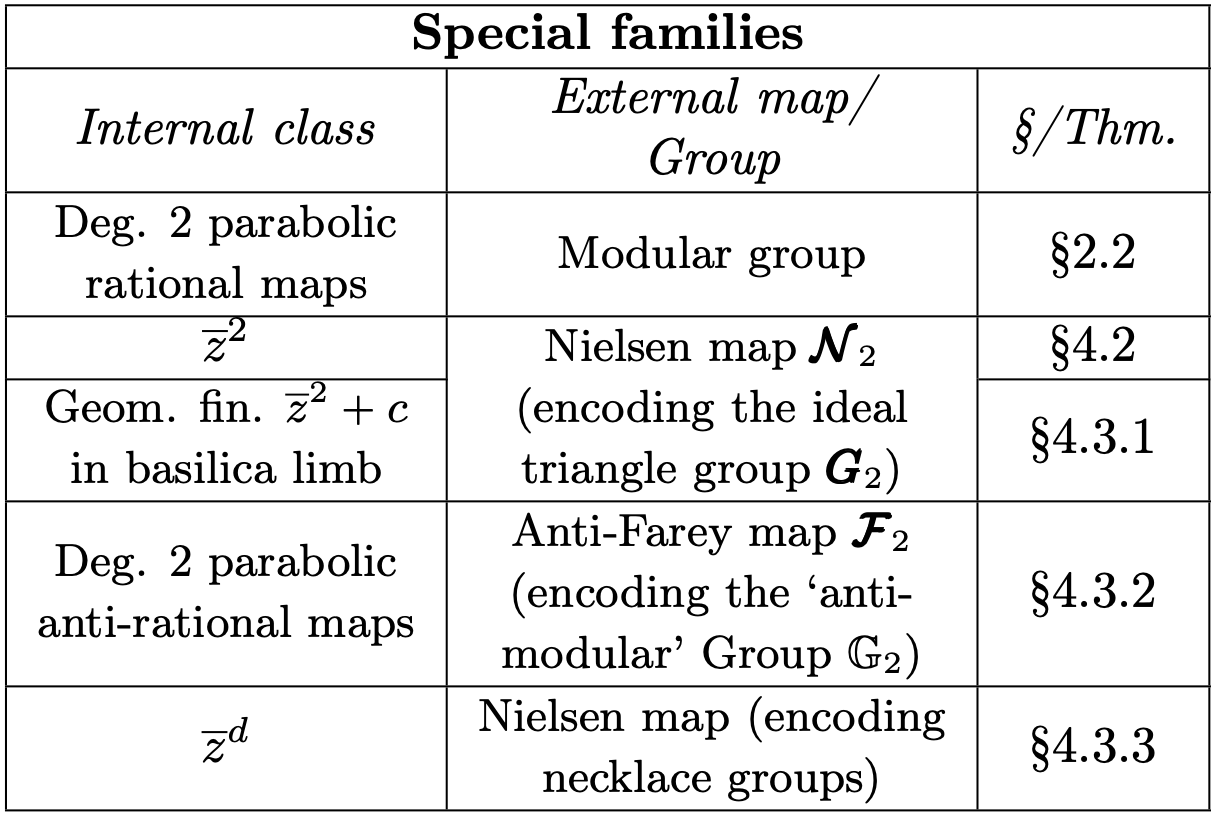}

\ \includegraphics[width=0.8\linewidth]{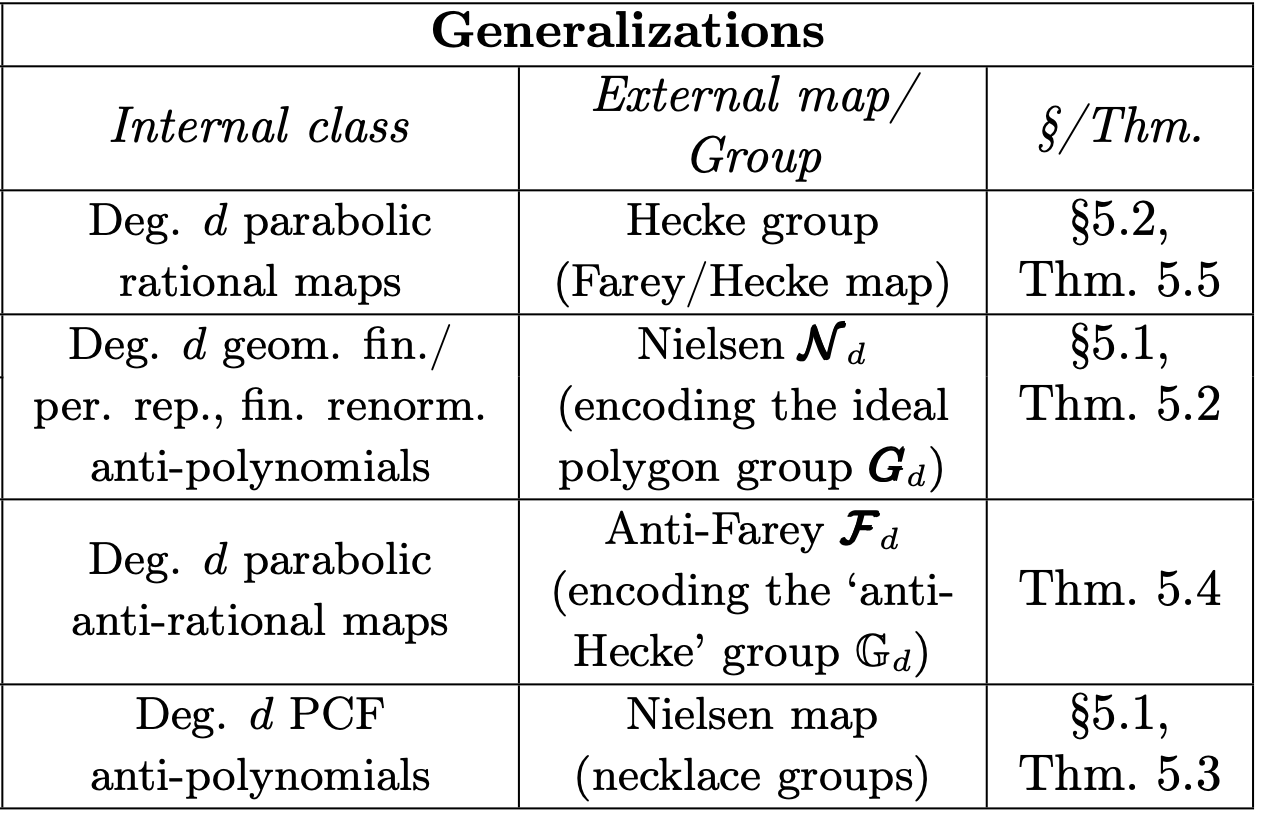}
\end{center}

\section{Mating phenomena in the antiholomorphic world: Schwarz reflection dynamics.}\label{schwarz_sec}

A new link between rational dynamics and Kleinian groups emerged in the late 2010s through the investigation of certain problems having their roots in statistical physics and complex analysis. In an influential work by Lee and Makarov \cite{LM16}, the study of equilibrium states of certain $2$-dimensional Coulomb gas models was connected to the iteration of a new class of anti-holomorphic maps, namely \emph{Schwarz reflection maps} associated with \emph{quadrature
domains} \cite{AS76}. 
A systematic exploration of Schwarz dynamical systems was launched in \cite{LLMM23a,LLMM25,LLMM21}, where it transpired that they often arise as \emph{conformal matings} of Kleinian reflection groups and anti-holomorphic rational maps (anti-rational maps in short), giving rise to a framework for such combinations. The underlying combination principle also had a significant influence on the holomorphic side of the mating story, allowing for extensions of some of the results mentioned in \S\ref{BP_sec} to wider classes of Fuchsian groups uniformizing genus zero orbifolds (see \S\ref{holo_gen_mating_subsec}).

\subsection{Quadrature domains and Schwarz reflection maps.}\label{schwarz_basic_prop_subsec} 
A proper sub-domain $\cU$ of the Riemann sphere $\widehat{\mathbb{C}}$ is called a \emph{quadrature domain} if it admits a continuous function $\sigma:\overline{\cU}\to\widehat{\mathbb{C}}$ such that $\sigma$ is anti-meromorphic on $\cU$ and $\sigma=\mathrm{id}$ on the boundary $\partial\cU$.
The map $\sigma$ is called the \emph{Schwarz reflection map} of $\cU$. We refer the reader to \cite{LM16} for various equivalent definitions of quadrature domains in terms of area/quadrature identities and Cauchy transforms.
The simplest example of a quadrature domain is given by a round disk, whose Schwarz reflection map is the standard anti-M{\"o}bius reflection in the boundary circle. 
\smallskip

\noindent\textbf{An important characterization.} Simply connected quadrature domains admit a simple description.
A simply connected domain $\cU\subsetneq\widehat{\C}$ (with $\infty\notin\partial\cU$ and $\mathrm{int}(\overline{\cU})=\cU$) is a quadrature domain iff the Riemann uniformization $R:\mathbb{D}\to\cU$ extends to a rational function of $\widehat{\C}$ \cite[Theorem~1]{AS76}. 
In this case, the Schwarz reflection $\sigma$ of $\cU$ is given by $R\circ\eta^-\circ(R\vert_{\D})^{-1}$, where $\eta^-(z):=1/\overline{z}$. 
Moreover, if $\deg{(R)}=d$, then $\sigma:\sigma^{-1}(\cU)\to\cU$ is a (branched) covering of degree $(d-1)$, and $\sigma:\sigma^{-1}(\Int{\cU^c})\to\Int{\cU}^c$ is a (branched) covering of degree~$d$.
\smallskip

\noindent\textbf{Piecewise Schwarz reflections, and a dynamical partition.}
Consider a \emph{piecewise Schwarz reflection map}
$\sigma: \bigcup_{i=1}^k \overline{\cU_i}\to\widehat{\C},\ z\mapsto \sigma_i(z),\ \textrm{for}\ z\in\overline{\cU}_i$; 
where $\{\cU_i\}$ are pairwise disjoint quadrature domains, and $\{\sigma_i\}$ are the associated Schwarz reflection maps. The union $\cU:=\bigcup_{i=1}^k \cU_i$ is called a \emph{quadrature multi-domain}. The singular points on $\partial\cU$ are cusps and double points (cf. \cite{Sak91,Gus88}), we denote this singular set by $\mathscr{S}$. 
We define the \emph{droplet}
$T:=\widehat{\C}\setminus \cU$, and the \emph{desingularized droplet} $T^0:= T\ \setminus\mathscr{S}$, see Figures~\ref{deltoid_fig},~\ref{necklace_mating_fig}. The sphere $\widehat{\C}$ admits a natural $\sigma$-invariant partition into the \emph{non-escaping/filled Julia} set $K(\sigma)$, the set of points that never leave $\cU\sqcup\mathscr{S}$, and the \emph{escaping set} $T^\infty(\sigma):=\widehat{\C}\setminus K(\sigma)$, consisting of points that eventually escape to the desingularized droplet $T^0$ (cf. \cite[\S 3.3]{LM23}). The common boundary of $T^\infty(\sigma)$ and $K(\sigma)$ is called the \emph{limit set}, and is denoted by $\Lambda(\sigma)$. The conformal class of $\sigma\vert_{K(\sigma)}$ (respectively, the conformal model of $\sigma$ on $T^\infty(\sigma)$) will be referred to as the \emph{internal class} (respectively, \emph{external map}) of $\sigma$.
We will now explicate a few examples, where the internal class of $\sigma$ is given by an anti-rational map, while the external map arises from a Kleinian reflection group. 
\subsection{The deltoid example: mating $\overline{z}^2$ with the ideal triangle reflection group $\pmb{G}_2$.}\label{deltoid_subsec}
The mating phenomenon displayed by Schwarz reflection maps was first instantiated in \cite{LLMM23a} in the case of the so-called deltoid, which we now describe briefly (see \cite[\S 4.1]{LM23} for a detailed account). This rather simple example, which can be regarded as an object halfway between quasi-Fuchsian groups \cite{Ber60} and quasi-Blaschke products \cite[\S 7.4]{BF14}, serves as a model for many of the combination results covered in this survey.

The rational map $R(z)=z+1/(2z^2)$ is injective on $\overline{\mathbb{D}}$, and hence $\cU=R(\mathbb{D})$ is a Jordan quadrature domain. Notably, the desingularized droplet $T^0$ takes the shape of an ideal hyperbolic triangle in this case, and the boundary of $T^0$ is the classical deltoid curve. 
The restriction $\sigma:\sigma^{-1}(\cU)\to\cU$ is a degree two proper antiholomorphic map between topological disks. However, the containment $\sigma^{-1}(\cU)\subset\cU$ is not compact; indeed, the domain $\sigma^{-1}(\cU)$ touches the boundary $\partial\cU$ at the deltoid cusps. This is an example of a pinched \emph{antiholomorphic polynomial-like map or a-p-l map} (see \S\ref{internal_external_sec}), and hence the dynamics of $\sigma$ splits into an external map and an internal class. 
\smallskip

\begin{figure}[htbp]
\captionsetup{width=0.96\linewidth}
\begin{tikzpicture}
\node at (0,0) {};
\node[anchor=south west,inner sep=0] at (0,0) {\includegraphics[width=0.96\textwidth]{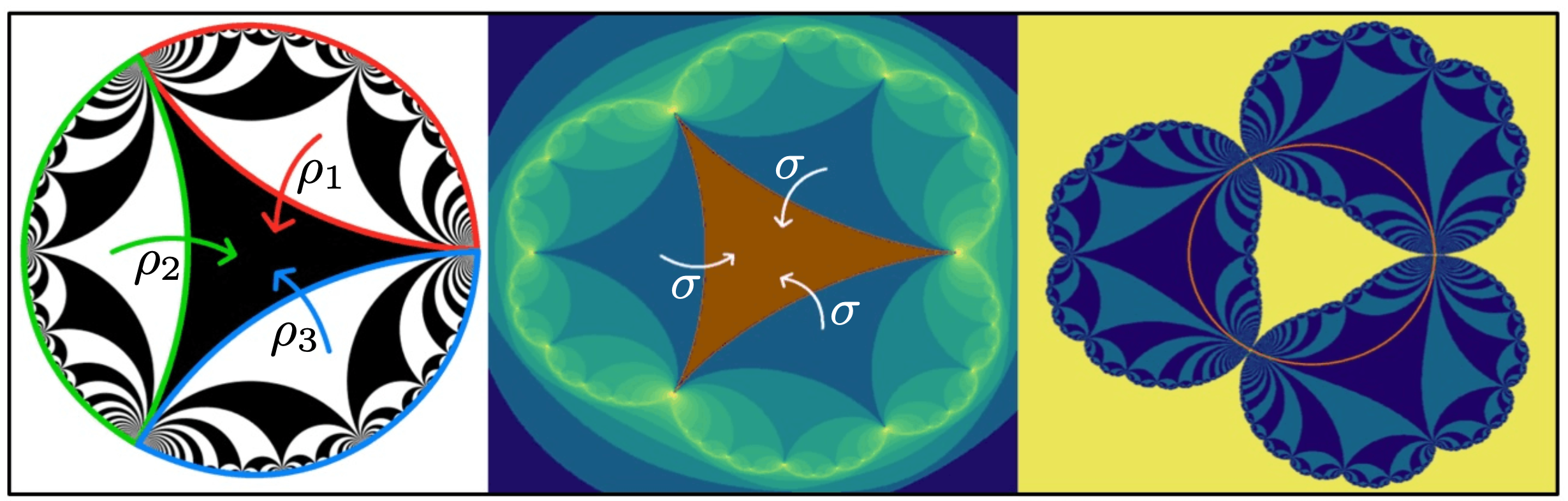}};
\node[anchor=south west,inner sep=0] at (1,-4.5) {\includegraphics[width=0.8\textwidth]{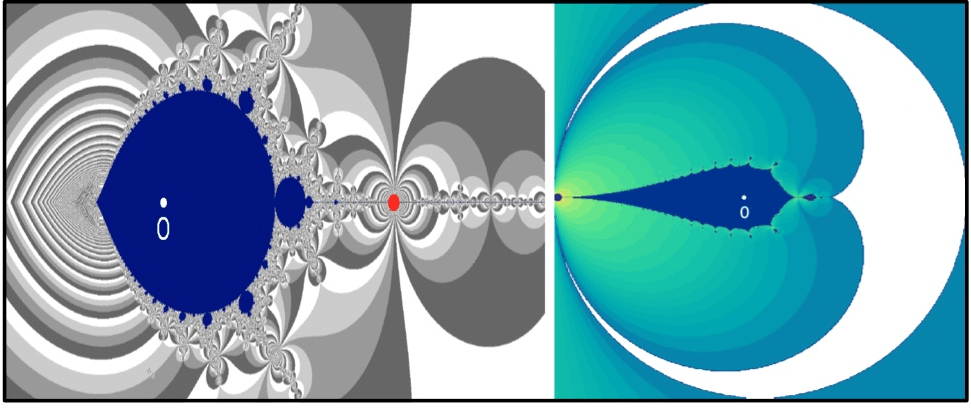}};
\end{tikzpicture}
  \caption{Top left: The Nielsen map $\pmb{\cN}_2$ of the ideal triangle group $\pmb{G}_2$. Top middle: The brown region is the droplet and its complement is the quadrature domain $\cU$ of \S\ref{deltoid_subsec}. The tiling set $T^\infty(\sigma)$, which resembles the tiling of $\D$ by ideal triangles, is the Jordan domain bounded by the green fractal. Top right: A branch of the associated correspondence $\mathfrak{C}$ preserves a yellow component and is conformally conjugate to $\overline{z}^2$. The branches of $\mathfrak{C}$ act as anti-conformal self-maps on the green/blue region, generating a properly discontinuous group action conformally equivalent to $\pmb{G}_2$. Bottom left: Part of the connectedness locus (and exterior tessellation) of the C\&C family of \S\ref{circle_and_cardioid_subsec}. Bottom right: Mating of $\overline{z}^2-1$ and $\pmb{\cN}_2$ in the C\&C family (droplet in white).}
  \label{deltoid_fig}
\end{figure}

\noindent\textbf{Nielsen external map.} The escaping set $T^\infty(\sigma)$ admits a tiling by the iterated $\sigma$-preimages of $T^0$. In fact, the action of $\sigma$ on $T^\infty(\sigma)$ is reminiscent of the \emph{ideal triangle reflection group} $\pmb{G}_2$, generated by the reflections in the sides of a (closed) hyperbolic triangle $\Pi$ in $\D$. Specifically, the dynamics of $\sigma$ on the tiling set is conformally conjugate to the \emph{Nielsen map} $\pmb{\cN}_2:\D\setminus\Int{\Pi}\to\D$ (associated with $\pmb{G}_2$), which acts on the three hyperbolic half-planes as reflections in the corresponding sides of $\Pi$ (see Figure~\ref{deltoid_fig}). In other words, the external map of the degenerate p-l restriction of $\sigma$ is given by $\pmb{\cN}_2$. The Nielsen map $\pmb{\cN}_2$ extends to the circle $\mathbb{S}^1$, is \emph{orbit equivalent} to $\pmb{G}_2$, and is an antiholomorphic analog of the so-called \emph{Bowen-Series maps} \cite{BS79}. 
\smallskip

\noindent\textbf{Internal class.} Classical straightening theorems do not apply to degenerate a-p-l maps.
In fact, the map $\sigma$ displays parabolic behavior near the cusps of the deltoid (which is to be expected from the fact that the group $\pmb{G}_2$ has parabolic elements); and hence, the dynamics of $\sigma$ around its non-escaping set is \emph{not} hybrid equivalent to the dynamics of an anti-polynomial. In \cite{LLMM23a}, tools from classical Fatou-Julia theory of rational maps (combined with puzzle techniques) were employed to show that the internal class of $\sigma$ is given by the map $\overline{z}^2\vert_{\overline{\D}}$. Thus, the dynamics of $\sigma$ combines the map $\pmb{\cN}_2$ with the anti-polynomial $\overline{z}^2$.  
A key fact that allows for this combination is that $\pmb{\cN}_2$ is an expansive covering of $\mathbb{S}^1$ that is topologically conjugate to $\overline{z}^2\vert_{\mathbb{S}^1}$. This compatibility between $\overline{z}^2$ and the Nielsen map $\pmb{\cN}_2$ acts as a bridge between reflection groups and anti-polynomial dynamics. 
\smallskip

\noindent\textbf{Parabolic vs. hyperbolic  mating, and non-quasisymmetric welding.} It is worth pointing out that existence of parabolic dynamics and degenerate a-p-l behavior are essential characteristics of the Schwarz reflection maps appearing in our mating framework. These features, which we already encountered in \S\ref{BP_sec}, will be a recurring theme in this survey.
Let us also mention that the limit set $\Lambda(\sigma)$ of the deltoid Schwarz reflection is a cuspidal Jordan curve, and hence not a quasicircle. The welding homeomorphism of $\Lambda(\sigma)$ is a conjugacy between the expanding map $\overline{z}^2$ and the parabolic map $\pmb{\cN}_2$, and is thus non-quasisymmetric \cite[\S 5]{LMMN25}.
\smallskip

\noindent\textbf{Antiholomorphic correspondence from deltoid reflection.} The Schwarz reflection $\sigma$ of the deltoid can be lifted, via the rational map $R$, to define an antiholomorphic correspondence $\mathfrak{C}$ on the sphere. Roughly speaking, we define a multi-valued map on $\widehat{\C}$, whose local branches are antiholomorphic functions, such that $z\mapsto w$ if the Schwarz reflection $\sigma$ or one of its backward branches carries $R(z)$ to $R(w)$. Concretely, the correspondence $\mathfrak{C}$ takes the form $\mathrm{Cov}_0^R\circ\eta^-$. 
The dynamical plane of $\mathfrak{C}$ is depicted in Figure~\ref{deltoid_fig}. The description of $\sigma$ as a mating between $\pmb{\cN}_2$ and $\overline{z}^2$ can be transferred to the dynamical plane of the correspondence via the map $R$. This shows that, with an appropriate notion of mating extending Definition~\ref{modular_mating_def}, the correspondence $\mathfrak{C}$ itself is a mating of the $\pmb{G}_2$ and $\overline{z}^2$ (see Definition~\ref{gen_mat_corr_def}; cf. \cite[Appendix~B]{LLMM21}, \cite[Theorem~1.9]{LLM24}).

\subsection{Special families of Schwarz reflections.}\label{special_fam_schwarz_subsec}
 In the early stage of development of Schwarz reflection dynamics, a number of families, generated by explicit univalent restrictions of rational maps, were studied. 
 We collect some key features of these families below (see \cite{LLMM21,LLMM23a,LLMM25} for proofs, or \cite{LM23} for a more thorough treatment).

\subsubsection{The Circle-and-Cardioid family: mating quadratic anti-polynomials with $\pmb{G}_2$.}\label{circle_and_cardioid_subsec}

We now come to the simplest family of quadratic Schwarz reflections; namely, the \emph{Circle-and-Cardioid (or C\&C) family}. These are piecewise Schwarz reflections $\sigma_a$ in the quadrature multi-domains $\cU_a$, where $\cU_a$ is the union of the main hyperbolic component $\heartsuit$ of the Mandelbrot set (i.e., $\heartsuit=R_1(\D)$, where $R_1(z)=z/2-z^2/4$), and the exterior of a round disk $B(a,r_a)$ such that $\heartsuit\subset B(a,r_a)$ and $\partial B(a,r_a)\cap\partial\heartsuit$ is a singleton (see Figure~\ref{deltoid_fig}). As the parameter $a$ runs over a slit-plane, we obtain a one-parameter family of Schwarz reflections $\{\sigma_a\}$. 
\smallskip

\noindent\textbf{Degenerate anti-polynomial-like maps.} As in the deltoid case, the restrictions $\sigma_a:\sigma_a^{-1}(\cU_a)\to\cU_a$ are quadratic degenerate a-p-l maps (with parabolic dynamics near the cardioid cusp and the point of tangency of $\partial\heartsuit$ and $\partial B(a,r_a)$).
Each such degenerate a-p-l map has a unique critical point (coming from the Schwarz reflection of the cardioid). The \emph{connectedness locus} of the family consists of maps with connected non-escaping sets $K(\sigma_a)$, or equivalently, of maps for which the unique critical point does not escape to the desingularized droplet $T^0$ under iteration. The complement of the connectedness locus is called the \emph{escape locus}.

\noindent\textbf{Internal class, external map, and matings.} The Tricorn is the connectedness locus of quadratic anti-polynomials $\{p_c(z)=\overline{z}^2+c\}_{c\in\C}$  (see \cite{HS14,IM16,MNS17,IM21} for background and structural differences between the Mandelbrot set and the Tricorn). In \cite{LLMM25}, a Fatou-Julia theory was developed for maps in the C\&C family, and it was proved that for each geometrically finite $p_c$ in the \emph{basilica limb} of the Tricorn, there exists a unique $\sigma_a$ in the connectedness locus of the C\&C family admitting $p_c$ as its internal class and $\pmb{\cN}_2$ (cf. \S\ref{deltoid_subsec}) as its external map.
\smallskip

\noindent\textbf{Tessellation of escape locus and combinatorial equivalence of parameter spaces.} Using a phase-parameter duality (as in the quadratic polynomial case), it was shown in \cite{LLMM25} that the location of the escaping critical value in appropriate uniformizing coordinates yields a dynamically natural uniformization of the escape locus by the disk. Thus, the escape locus admits a tessellation conformally equivalent to the tiling of $\D$ induced by the group $\pmb{G}_2$ (see Figure~\ref{deltoid_fig}). This external structure was exploited to establish that a combinatorial model of the connectedness locus of the C\&C family is homeomorphic to such a model of the basilica limb of the Tricorn. 

These results were extended to cover generic geometrically infinite anti-polynomials in \cite{LLM24}. There, these single-valued matings were also promoted to antiholomorphic correspondences on nodal spheres that capture the actions of the anti-polynomials as well as the full structure of the group $\pmb{G}_2$ (see \S\ref{antiholo_gen_mating_subsec}).

\subsubsection{Antiholomorphic counterpart of the family $\cF_a$.}\label{cubic_chebyshev_subsec}

Our next family of quadratic Schwarz reflections is motivated by the family $\cF_a$ treated in \S\ref{bpl_family_subsec}. Note that the univalent restriction of a rational map on any round disk yields a quadrature domain by the mechanism in \S\ref{schwarz_basic_prop_subsec}. 
Let $Q(z)=z^3-3z$ be the cubic Chebyshev polynomial defined in \S\ref{bpl_family_subsec}. Let $\lambda$ be such that the set $\mathscr{B}_\lambda:=\{z:\vert z-\lambda\vert\leq\vert 1-\lambda\vert\}$, having the marked critical point $1$ on their boundary, is a  domain of univalency for $Q$. Then $\cU_\lambda:=Q(\Int{\mathscr{B}_\lambda})$ is a quadrature domain. This gives rise to a family of quadratic Schwarz reflections $\sigma_\lambda:\overline{\cU}_\lambda\to\widehat{\C}$, which we call the \emph{cubic Chebyshev family} of Schwarz reflections. The restriction $\sigma_\lambda:\sigma_\lambda^{-1}(\cU_\lambda)\to\cU_\lambda$ is a quadratic degenerate a-p-l map having a (free) simple critical point: the non-escaping set $K(\sigma_\lambda)$ is connected (equivalently, $T^\infty(\sigma_\lambda)$ is a topological disk) precisely when this critical point does not escape to the desingularized droplet $T^0(\sigma_\lambda)$.
\smallskip

\noindent\textbf{Internal class and hybrid conjugacy.}
The degenerate a-p-l maps arising from the cubic Chebyshev family $\{\sigma_\lambda\}$ have a unique parabolic fixed point (unlike the deltoid or C\&C case). This suggests that such maps can be straightened by means of a quasiconformal surgery to parabolic anti-rational maps with a persistently parabolic fixed point. In \cite{LLMM21}, such a surgery was designed using a classical result of Warschawski (cf. \cite{War42}) to demonstrate that the said degenerate a-p-l restrictions are hybrid conjugate to certain quadratic parabolic anti-rational maps.
\smallskip

\noindent\textbf{Anti-Farey external map.}
Each Schwarz reflection in our family has an escaping critical point; i.e., a (passive) critical point that escapes to $T^0(\sigma_\lambda)$ in one step. While this seems to obstruct a group structure in the tiling set, the external maps of the associated degenerate a-p-l restrictions still arise from certain groups. For connected $K(\sigma_\lambda)$, the dynamics of $\sigma_\lambda$ on $T^\infty(\sigma_\lambda)$ is modeled by a factor of the Nielsen map $\pmb{\cN}_2$. Since the map $\pmb{\cN}_2$ commutes with $M_\omega(z)=\omega z$, $\omega=e^{2\pi i/3}$, it descends to a factor map on the quotient $\D/\langle M_\omega\rangle$. The resulting factor map, the \emph{degree $2$ anti-Farey map} $\pmb{\cF}_2$, has a unique parabolic fixed point on $\mathbb{S}^1$ and an escaping critical point. The map $\pmb{\cF}_2$ serves as the external map of the degenerate a-p-l maps coming from the connectedness locus of the cubic Chebyshev family $\sigma_\lambda$ (cf. \cite{LLMM21}, see Figure~\ref{necklace_mating_fig}). The existence of hybrid conjugacies between the maps $\sigma_\lambda$ and parabolic anti-rational maps can be attributed to the fact that $\pmb{\cF}_2$ has a unique parabolic fixed point, and hence is quasisymmetrically conjugate to a parabolic anti-Blaschke product (cf. \S\ref{gen_matings_sec}).
\smallskip

\noindent\textbf{Correspondences as matings.}
As in \S\ref{deltoid_subsec}, one can lift the maps $\sigma_\lambda$ to define antiholomorphic correspondences $\mathrm{Cov}_0^Q\circ\eta^-_\lambda$, where $\eta^-_\lambda$ is the reflection in the circle $\partial\mathscr{B}_\lambda$.
Note that this is the antiholomorphic version of the family $\cF_a=J_a \circ \mathrm{Cov}_0^Q$ (see \S\ref{bpl_family_subsec}), with the conformal involution $J_a$ replaced by the anti-conformal involution $\eta^-_\lambda$. For the analogies, see \cite[Remark~1]{BL22}, and compare \cite[Figure 8]{BL20} with the right sides of \cite[Figures 13,16]{LLMM21}.
The mating structure of the Schwarz reflections in the connectedness locus can be pulled back by $Q$ to the correspondence planes yielding a description of the correspondences as matings between quadratic parabolic anti-rational maps and the group $\mathbb{G}_2=\langle\pmb{G}_2,M_\omega\rangle$ (see Figure~\ref{necklace_mating_fig} for the group structure of the correspondence, and  \cite{LLMM21,LMM24} for more details). Note that $\mathbb{G}_2\cong C_2*C_3$, and can be regarded as an anti-conformal analog of $\PSL_2(\Z)$.
\smallskip

\noindent\textbf{Parameter space results.} An analysis of the regularity properties of the straightening map and a dynamically natural tessellation of the escape locus (given by the conformal position of the escaping free critical value) revealed a combinatorial equivalence between the connectedness locus of this family and the parabolic Tricorn.
\subsubsection{The family $\Sigma_d^\ast$: mating $\overline{z}^d$ with necklace groups.}\label{sigma_d_subsec}

We now talk about a class of Schwarz reflections that combine the power map $\overline{z}^d$, $d\geq 2$, with reflection groups whose limit sets are not Jordan curves. 
\smallskip

\noindent\textbf{A Bers slice of reflection groups.} The \emph{regular ideal $(d+1)$-gon reflection group $\pmb{G}_d$} is the (discrete) group generated by reflections in the sides of a regular ideal $(d+1)$-gon in the hyperbolic plane. The notion of a Nielsen map naturally extends to this case (see \S\ref{deltoid_subsec}), yielding a piecewise anti-M{\"o}bius map $\pmb{\cN}_d$ (defined on a pinched neighborhood of $\mathbb{S}^1$), that is orbit equivalent to $\pmb{G}_d$, and topologically (but non-quasisymmetrically) conjugate to $\overline{z}^d\vert_{\mathbb{S}^1}$. For $d\geq 3$, the group $\pmb{G}_d$ admits quasiconformal deformations, by moving the ideal vertices. The \emph{Bers slice} $\mathfrak{B}(\pmb{G}_d)$ consists of quasiconformal deformations of $\pmb{G}_d$ supported on $\D$ (see \cite[\S 2.2]{LMM22}). A \emph{necklace group} is a reflection group lying in the closure of $\mathfrak{B}(\pmb{G}_d)$, where the closure is taken in the character variety. Equivalently, a necklace group is the group generated by reflections in the circles of a finite circle packing with a $2$-connected, outerplanar contact graph \cite[\S 3.4]{LLM22a}. The Nielsen map $\cN_G$ of a necklace group $G$ is defined on the union of the $d+1$ (closed) disks of the underlying circle packing as reflections in the corresponding circles. The limit set $\Lambda_G$ is the quotient of $\mathbb{S}^1$ under an $\pmb{\cN}_d$-invariant geodesic lamination, and hence, there is a continuous surjection $\Phi_G:\mathbb{S}^1\to\Lambda_G$ (called the \emph{Cannon-Thurston} map, cf. \cite{MS13,Mj18}), that semi-conjugates $\pmb{G}_d$ to~$\cN_G$. 
The existence of $\Phi_G$ makes $\cN_G\vert_{\Lambda_G}$ dynamically compatible with $\overline{z}^d\vert_{\mathbb{S}^1}$.
\begin{figure}[htbp]
\captionsetup{width=0.96\linewidth}
  \centering
  \includegraphics[width=0.96\linewidth]{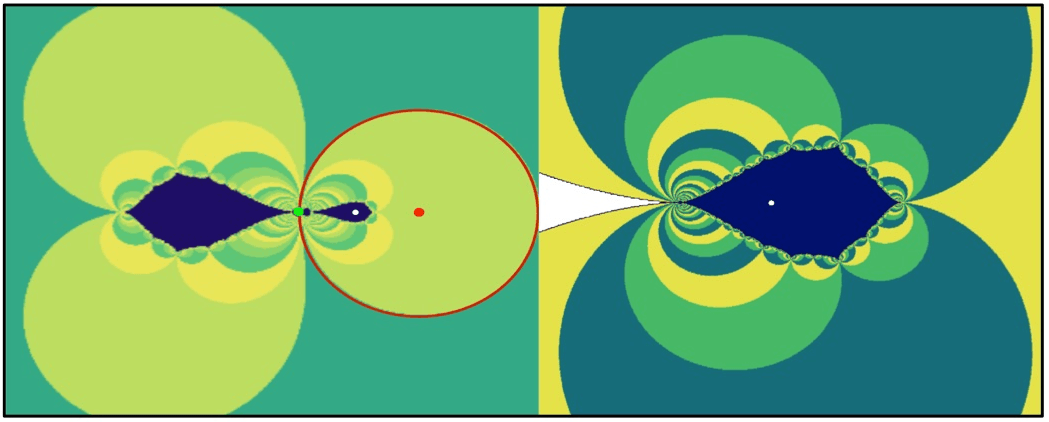}
    \includegraphics[width=0.8\linewidth]{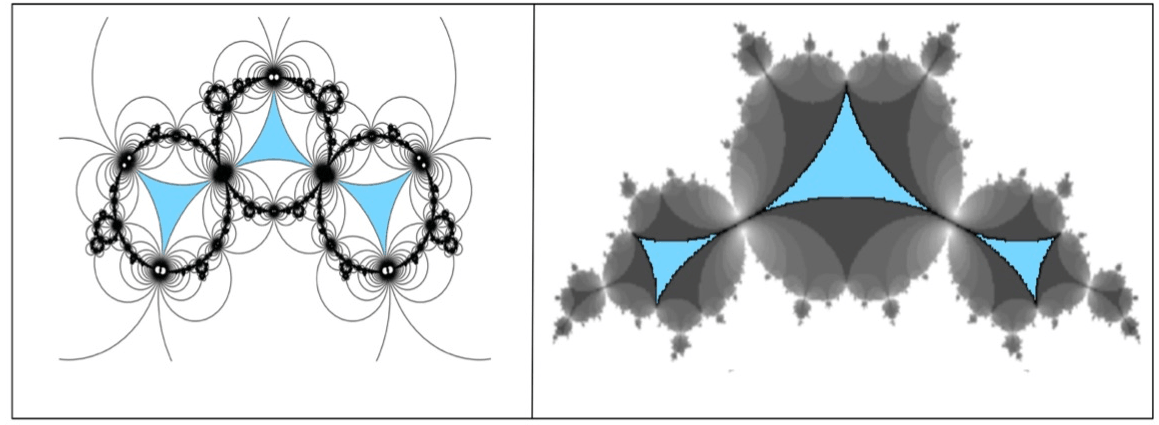};
  \caption{Top: The dynamical planes of a Schwarz reflection in the cubic Chebyshev family $\{\sigma_\lambda\}$ (right) and the associated correspondence (left), as in \ref{cubic_chebyshev_subsec}. The red circle bounds a disk where $Q$ is injective. Bottom: The limit sets of a necklace group $G$ generated by reflections in $5$ circles (left) and the Schwarz reflection that realizes the mating between $\cN_G$ and the anti-polynomial $\overline{z}^4$, as in \ref{sigma_d_subsec} (right). The blue region in the rightmost picture is the droplet.}
  \label{necklace_mating_fig}
\end{figure}
\smallskip

\noindent\textbf{Bers slices and $\Sigma_d^*$.} The space $\Sigma$ consists of normalized conformal maps $f:\overline{\D}^\complement\to\widehat{\C}$ with a simple pole at $\infty$ (cf. \cite{Dur83}). The subcollection 
$$
\Sigma_d^* := \left\{ f(z)= z+a_1/z + \cdots +a_d/z^d : a_d=-1/d\textrm{ and } f|_{\mathbb{D}^*} \textrm{ is conformal}\right\}
$$ 
is a compact subspace of $\Sigma$. For a rational map $R\in\Sigma_d^*$, the complement of the quadrature domain $\cU=R(\overline{\D}^\complement)$ is a `tree of polygons'. The topology of this tree of polygons can be matched by an appropriate circle packing giving rise to a necklace group $G$ (see Figure~\ref{necklace_mating_fig}). The associated Schwarz reflection $\sigma\equiv R\circ\eta^-\circ(R\vert_{\D^\complement})^{-1}:\overline{\cU}\to\widehat{\C}$ has an order $d$ pole at $\infty$, and hence combines the power map $\overline{z}^d$ with the Nielsen map $\cN_G$. (In the rightmost picture of Figure~\ref{necklace_mating_fig}, the white region is where $\sigma$ is conformally conjugate to $\overline{z}^4\vert_\D$, while the action of $\sigma$ on the gray/black region is conformally conjugate to $\cN_G$ for a necklace group $G$ whose limit set is displayed next to it.)

Using a combination of tools from polynomial dynamics, reflection groups, and extremal length arguments, it was shown in \cite{LMM22} that matings of all groups in $\overline{\mathfrak{B}(\pmb{G}_d)}$ with $\overline{z}^d$ are realized by Schwarz reflections arising from $\Sigma_d^*$, and moreover, this mating operation provides a natural homeomorphism between the spaces $\overline{\mathfrak{B}(\pmb{G}_d)}$ and $\Sigma_d^*$.

\section{From special families to a combination program.}\label{gen_matings_sec}


Most of the results discussed hitherto can be regarded as \emph{unmating/decomposition theorems}, where one breaks a given conformal dynamical system into simpler pieces. We now reverse the point of view, and address the question of mateability of  groups $G$ and maps $R$. In general, this seems to be a hard problem due to the intrinsic mismatch between group dynamics (invertible, many generators) and map dynamics (non-invertible, single generator). However, several examples mentioned so far suggest the following four-step program to produce matings of maps and groups realized as correspondences. This strategy lies at the heart of many of the general results that we will mention in this section.

\noindent \textbf{I.} Given a group $G$, come up with a continuous, piecewise-analytic map $A_G$, called a \emph{virtually mateable map}, that remembers the essential dynamics of $G$, and is compatible with rational dynamics in an appropriate way.

\noindent \textbf{II.} Use surgery and 'mating continuity' techniques to mate $A_G$ with $R$, producing a single-valued map $F$ that combines $R$ and $G$ in a certain weak sense. The map $F$ can often be interpreted as a degenerate (anti-)polynomial-like map, whose \emph{internal class} is given by $R$, while its \emph{external map} is $A_G$ (cf. \cite{BLLM25,LLM24}). 

\noindent \textbf{III.} Identify the analytic map $F$ as an algebraic function.

\noindent \textbf{IV.} Use the above algebraic description to promote $F$ to an algebraic correspondence on a compact surface.\\
The following definition extends Definition~\ref{modular_mating_def} to incorporate groups with multiple parabolics.
\begin{definition}\label{gen_mat_corr_def}
Let $P$ be a degree $d$ (anti-)polynomial with connected Julia set and let $G$ be a discrete subgroup of (anti-)conformal automorphism group of the unit disk $\D$.
Let $\mathfrak{C}$ be an (anti-)holomorphic correspondence on a (possibly noded) sphere $\cS$.
We say $\mathfrak{C}$ a {\em mating} of $P$ and $G$ if there is a $\mathfrak{C}$-invariant partition 
$\cS=\cT\sqcup\cK$ such that the following hold.

1. On $\cT$, the dynamics of $\mathfrak{C}$ is equivalent to the action of a group of (anti-)conformal automorphisms acting properly discontinuously. Further, $\cT/\mathfrak{C}$ is biholomorphic to $\D/G$.
	
2. $\cK=\cK_-\cup\cK_+$, each $\cK_\pm$ is a copy of the filled Julia set $K(P)$ of $P$, such that $\cK_-,\cK_+$ intersect in (at most) finitely many points. Furthermore, $\mathfrak{C}$ has a forward (respectively, backward) branch carrying $\cK_-$ (respectively, $\cK_+$) onto itself with degree $d$, and this branch is conformally (respectively, anti-conformally) conjugate to $P\vert_{K(P)}$. 
\end{definition}

\subsection{General matings in the antiholomorphic world.}\label{antiholo_gen_mating_subsec}
We now discuss general mating results between anti-polynomials and reflection (or anti-Hecke) groups in arbitrary degree (see \cite{LM23} for a detailed exposition).

\noindent\textbf{Mating anti-polynomials with Fuchsian reflection groups.} Generalizing the setup of \S\ref{cubic_chebyshev_subsec}, we define $\mathbb{G}_d:=\langle\pmb{G}_d,M_\omega\rangle$, where $\omega=e^{2\pi i/d+1}$, and call it the \emph{anti-Hecke group}.
\begin{theorem}\upshape\cite[Theorem~1.9]{LLM24}\label{polygon_reflection_grp_anti_poly_corr_thm}
Let $G\in\{\pmb{G}_d,\mathbb{G}_d$\}.
Let $P$ be a degree $d$ anti-polynomial with connected Julia set which is either geometrically finite; or periodically repelling, finitely renormalizable.
Then there exists an antiholomorphic correspondence $\mathfrak{C}$ on a (nodal) sphere $\cS$ that is a mating of $P$ and $G$.
The correspondence $\mathfrak{C}$ is of the form $\mathrm{Cov}_0^R\circ \widehat{\eta}^-$, where $R:\cS\to\widehat{\C}$ is a branched cover and $\widehat{\eta}^-$ is an anti-conformal involution on $\cS$.
\end{theorem}
\noindent \emph{Key ideas of the proof.} The Nielsen map $\pmb{\cN}_d$ of the ideal $(d+1)$-gon reflection group $\pmb{G}_d$ is topologically conjugate to $\overline{z}^d$ on $\mathbb{S}^1$. As in the degree two case (see \S\ref{cubic_chebyshev_subsec}), the map $\pmb{\cN}_d$ commutes with multiplication by $\omega=e^{2\pi i/d+1}$, and hence descends to the degree $d$ \emph{anti-Farey} map $\pmb{\cF}_d$ on the quotient $\overline{\D}/\langle M_\omega\rangle$, where $M_\omega(z)=\omega z$. The map $\pmb{\cF}_d$ is also conjugate to $\overline{z}^d$ on $\mathbb{S}^1$. On the other hand, if a degree $d$ anti-polynomial $P$ has a locally connected Julia set, then $P\vert_{J(P)}$ is a factor of $\overline{z}^d\vert_{\mathbb{S}^1}$.
The above observations allow one to construct a \emph{topological} degenerate a-p-l map having $P$ as its internal class, and $\pmb{\cN}_d$ or $\pmb{\cF}_d$ as its external map. As the conjugacy between $\overline{z}^d$ and $\pmb{\cN}_d$ (or $\pmb{\cF}_d$) is not quasisymmetric, turning the topological mating into a conformal one is challenging. 
To address this problem, it was proved in \cite{LMMN25} (cf. \cite{LN24a}) that the circle homeomorphism conjugating $\overline{z}^d$ to $\pmb{\cN}_d$ (or $\pmb{\cF}_d$) extends continuously to $\D$ as a \emph{David homeomorphism} (or homeomorphisms of exponentially integrable distortion, see \cite{Dav88,AIM09}). This extension result and a surgery argument involving David maps can be used to upgrade the above topological matings to conformal ones, under the additional assumption that $p$ is geometrically finite or subhyperbolic, see \cite[\S 7]{LMMN25}, \cite[\S 5]{LLM24} (this assumption is required to apply the \emph{David integrability theorem}). Further, such conformal matings turn out to be Schwarz reflections admitting degenerate a-p-l restrictions with $\pmb{\cN}_d$ (or $\pmb{\cF}_d$) as their external map. These spaces of degenerate a-p-l maps are compact \cite{LLM24,LMM24}, and a combination of puzzle and combinatorial continuity/rigidity techniques were employed in \cite{LLM24} to mate periodically repelling, finitely renormalizable anti-polynomials in the connectedness locus with the maps $\pmb{\cN}_d$, $\pmb{\cF}_d$. Finally, the algebraic description of Schwarz reflections (see \S\ref{schwarz_basic_prop_subsec}) permits one to lift these matings to construct the desired antiholomorphic correspondences (realizing matings between $P$ and $\pmb{G}_d, \mathbb{G}_d$) defined on possibly nodal spheres.\hspace{5cm} $\blacksquare$ 
\smallskip

\noindent\textbf{Mating anti-polynomials with necklace groups.} As mentioned in \S\ref{sigma_d_subsec}, the action of the Nielsen map $\cN_G$ of a general necklace group $G$ on $\Lambda_G$ is a factor of $\pmb{\cN}_d\vert_{\mathbb{S}^1}$. Hence, one can use the conjugacy between $\pmb{\cN}_d$ and $\overline{z}^d$ to topologically mate $\cN_G$ with anti-polynomials $P$ (with locally connected Julia set), as long as there is no \emph{Moore obstruction} (cf. \cite{MP12}).
Analogs of the following result are well-known for matings of two geometrically finite polynomials or two geometrically finite Bers boundary groups (cf. \cite[\S 9.3]{BF14}).
\begin{theorem}\upshape\cite{LMMN25}\label{necklace_pcf_anti_poly_mating_thm}
A degree $d$ pcf, hyperbolic anti-polynomial and a rank $d+1$ necklace group are conformally mateable iff they are topologically mateable, and the mating is realized as a Schwarz reflection.    
\end{theorem}
The proof uses the existence of equivariant homeomorphisms between limit sets of necklace groups and Julia sets of critically fixed anti-polynomials (see \S\ref{connections_sec}), tools from polynomial matings, and David surgery arguments. 
\smallskip

\noindent\textbf{Mating parabolic anti-rational maps with anti-Hecke groups.}
The anti-Farey map $\pmb{\cF}_d$ has a unique parabolic fixed point on $\mathbb{S}^1$ (cf. \S\ref{cubic_chebyshev_subsec}), and hence it is quasisymmetrically conjugate to a parabolic anti-Blaschke product. This fact enables one to generalize the family of \S\ref{cubic_chebyshev_subsec} to arbitrary degree. More precisely, one can apply a quasiconformal surgery to construct Schwarz reflections admitting degenerate a-p-l restrictions that have $\pmb{\cF}_d$ as the external map and suitable parabolic anti-rational maps as hybrid classes. The associated correspondences $\mathrm{Cov}_0^Q\circ\eta^-$ (where $Q$ is the uniformizing rational map for the Schwarz reflection) can be seen to be matings of parabolic anti-rational maps and the anti-Hecke group $\mathbb{G}_d:=\langle\pmb{G}_d,M_\omega\rangle$. 
\begin{theorem}\upshape\cite[Theorem~A]{LMM24}\label{anti_hecke_anti_poly_mating_thm}
Let $R$ be a degree $d$ anti-rational map with a parabolic fixed point of multiplier $1$ having a fully invariant and simply connected immediate basin of attraction. Then, there exists an antiholomorphic correspondence $\mathfrak{C}=\mathrm{Cov}_0^Q\circ\eta^-$ that is a mating of $\mathbb{G}_d$ and $R$ (where $Q$ is a degree $d+1$ polynomial that is univalent on $\overline{\D}$ and has a unique critical point on $\mathbb{S}^1$).
\end{theorem}
We refer the reader to \cite{LMM24,LMM25d}, \cite[\S 11]{LM23} for parameter space implications of the above mating result.

\subsection{General matings in the holomorphic world.}\label{holo_gen_mating_subsec}
The strategy outlined in the beginning of \S\ref{gen_matings_sec} can be implemented to mate several Fuchsian groups with  polynomials. 
\smallskip

\noindent\textbf{Mating parabolic rational maps with Hecke groups.}
The first examples of correspondences that are matings between Hecke groups and polynomials were constructed in \cite{BF1,BF2}. In \cite{BF1} it was conjectured that for each degree $d$ polynomial $P$ with a connected Julia set, there should exist a polynomial $Q$ of degree $d+1$ together with an involution $J$ such that the $d:d$ correspondence $J \circ Cov_0^Q$ realizes a mating of $P$ with the Hecke group $\mathcal H_{d+1}$. Implementing the strategy outlined in the beginning of this section,  this conjecture was confirmed in \cite{BLLM25}, after replacing `polynomials' with suitable `parabolic rational maps'.
Note that the principle of \S\ref{para_rat_mod_group_mating_subsec} is maintained: groups with a unique conjugacy class of parabolics admit quasiconformal matings with parabolic rational maps. 
\begin{theorem}\label{par_rat_hecke_thm}\upshape\cite{BLLM25}
Let $R$ be a degree $d$ rational map with a parabolic fixed point of multiplier $1$ having a fully invariant and simply connected immediate basin of attraction.
Then there exists a $d:d$ holomorphic correspondence $\cF$ on the Riemann sphere $\widehat{\C}$ which is a mating between $R$ and $\mathcal{H}_{d+1}$. Moreover, $\cF=J \circ Cov_0^Q$, where $J$ is a conformal involution and $Q$ is a degree $d+1$ polynomial.     
\end{theorem}
\smallskip

\noindent\textbf{Combining genus zero orbifolds with polynomials.}
Let $\mathfrak{F}$ be the class of finite-type hyperbolic orbifolds of genus zero with at most one order $2$, and at most one order $\nu\geq 3$ orbifold points.

\begin{theorem}\upshape\cite{LLM24,MM25}\label{genus_zero_orb_grp_poly_corr_thm}
Let $G$ be a Fuchsian group with $\Sigma=\D/G\in\mathfrak{F}$.
Let $d_\Sigma=1-2\alpha\cdot\chi_{\mathrm{orb}}(\Sigma)$; where $\alpha=\nu$ if $\Sigma$ has an order $\nu \geq 3$ orbifold point, and $\alpha=1$ otherwise.
Let $P$ be a degree $d_\Sigma$ polynomial with connected Julia set which is either geometrically finite; or periodically repelling, finitely renormalizable.
Then there exists a holomorphic correspondence $\mathfrak{C}=\mathrm{Cov}_0^R\circ \widehat{\eta}^+$ on a (possibly nodal) sphere $\cS$ that is a mating of $P$ and $G$ (where $R:\cS\to\widehat{\C}$ is a branched cover and $\widehat{\eta}^+$ is a conformal involution on $\cS$).
\end{theorem}

\noindent\emph{Technical distinctions from the antiholomorphic case.} While the core idea here is similar to the one used in Theorem~\ref{polygon_reflection_grp_anti_poly_corr_thm}, there are several additional hurdles. \textbf{I.} To construct virtually mateable maps for the group $G$, we resort to classical \emph{Bowen-Series maps}: these are piecewise M{\"o}bius circle maps cooked out of Fuchsian groups \cite{BS79}.
A Fuchsian punctured sphere group $G$ (possibly with an order $2$ orbifold point) admits a continuous Bowen-Series map that is topologically conjugate to some $z^d$ (see \cite{MM23a}, cf. \cite{Laz21}). This is no longer true when $G$ has torsion elements of order $\nu\geq 3$; in this case, one needs to work with a Bowen-Series map of a suitable index $\nu$ subgroup of $G$ (such that this map has a $\nu$-fold rotational symmetry), and then pass to $\nu$-fold factor map to get the desired virtually mateable map $A_G$. The map $A_G$ is called a \emph{factor Bowen-Series map} of $G$ (see \cite{MM25}, cf. \S\ref{antiholo_gen_mating_subsec}). 
\textbf{II.} A key feature of the conformal mating of a polynomial $P$ and a factor Bowen-Series map $A_G$ is that it induces an orientation-reversing involution on the piecewise-analytic boundary of its natural domain of definition, giving birth to a class of meromorphic maps called \emph{B-involutions}. Conformal welding arguments were employed in \cite[\S 14]{LLM24} (cf. \cite{MV25}) to prove algebraicity of B-involutions. Specifically, B-involutions admit an algebraic description similar to the one given in \S\ref{schwarz_basic_prop_subsec}, where the role of $\eta^-$ is now played by the conformal involution $\eta^+(z)=1/z$, and the role of $\D$ is played by a Jordan disk $\mathfrak{D}$ that is mapped inside out by $\eta^+$. (This gives a method to produce correspondences satisfying Bullett's combination criterion \cite[\S 3]{Bul00} in certain cases.) 
\textbf{III.} Compactness of the space of degenerate p-l maps having $A_G$ as their external map (which is necessary to go from geometrically finite to generic matings) is more subtle here, partly owing to lack of uniformity/regularity of the domain $\mathfrak{D}$. In \cite{LMM25a}, this was achieved using a delicate rescaling limit argument (involving the algebraic description of $B$-involutions). 
\smallskip

\noindent\textbf{Mating without parabolics.}
As mentioned in \S\ref{bpl_family_subsec}, matings between quadratic polynomials and representations in $\Int{\cD}$ (see \S\ref{bpl_family_subsec}) were constructed in \cite{BH00} using quasiconformal surgery.
This process has been
generalized to matings between degree $d$ polynomials and discrete representations of
$C_2*C_{d+1}$ with Cantor limit sets by Ratis Laude in \cite{Lau24}, producing correspondences which on the
complement of the forward and backward limit sets (a topological annulus)
have orbit space identical to that of a representation of $C_2*C_{d+1}$.
See also the thesis \cite{Lau25} for a mating construction between the fat basilica and a circle packing group living in~$\partial\cD$.
\smallskip

\noindent\textbf{Mating rational maps with $(p+1,q+1,\infty)$ triangle groups.}
These are realised by $pq:pq$ correspondences $\cF$ of the form
$Cov_
0^P*Cov_0^Q$ where $P, Q$ are polynomials of degree $p+1, q+1$ respectively, and which
have a common simple critical point. The construction generalises that of
correspondences $J
_a\circ Cov_0^Q$ discussed earlier (which is the case
$p=1$).
In the matings, the sphere is partitioned into an
invariant open set $\Omega$ and its complement, a one point union $
\Lambda_-\cup \Lambda_
+$. The set $\Lambda_
-$ is a copy of the filled Julia
set $K(g\circ f)$, and $\Lambda_
+$ is a copy of $K(f\circ g)$. The
correspondence $\cF$ has a branch $\Lambda_- \to \Lambda_
+$ conjugate to
$f:K(g\circ f) \to K(f\circ g)$, and $\cF^{-1}$ has a branch $\Lambda_
+ \to
\Lambda_
-$ conjugate to $g:K(f\circ g) \to K(g\circ f)$. The set $\Omega$ is a
copy of $\mathbb{H}$ carrying the parabolic representation of $C_{p+1}*C_{q+1}$
(which has fundamental domain a $(p+1,q+1,\infty)$-triangle union its reflection). See \cite{BL25}. 
\smallskip

\noindent\textbf{Correspondences on higher genus surfaces.}
A general combination theorem for finitely many Blaschke products and Fuchsian genus zero orbifold groups was presented in \cite{MV25}. In a nutshell, one starts with a rational map having several invariant Fatou components, and replaces its action on these invariant components with given Blaschke products and factor Bowen-Series maps (see \S\ref{holo_gen_mating_subsec}), thus producing a single-valued mating akin to a B-involution. The algebraic description of such matings involves a compact Riemann surface $\cS$ equipped with a hyperelliptic involution $\widecheck{\eta}$ (analog of $\eta^+(z)=1/z$) and a meromorphic map $\cR:\cS\rightarrow\widehat{\C}$. This leads to correspondences on surfaces of arbitrary genus realizing matings of Blaschke products and Fuchsian groups.
\begin{theorem}\upshape\cite[Theorem~A]{MV25}\label{corr_on_higher_genus_thm}
Given hyperbolic Blaschke products $\{B_j\}_{j=1}^r$, and Fuchsian groups $\{G_i\}_{i=1}^l$, with $\D/G_i\in \mathfrak{F}$, there exists an algebraic correspondence $\mathfrak{C}:=\mathrm{Cov}_0^{\cR}\circ\widecheck{\eta}\ $ on a hyperelliptic compact Riemann surface $\cS$ that combines the dynamics of $\{G_i\}_{i=1}^l$ and $\{B_j\}_{j=1}^r$ (where $\widecheck{\eta}$ and $\cR$ are as above).
\end{theorem}

\subsection{The inverse problem, Modular Multibrot sets, and Belyi polynomials.}\label{inverse_subsec}
Most correspondences realizing matings of rational maps and Kleinian groups can be written as $J\circ\mathrm{Cov}_0^Q$, for an involution $J$ and a rational map $Q$.
In the case of matings of Hecke groups with parabolic rational maps, the rational map $Q$ in question is a polynomial with a simple critical point that is fixed by the involution $J$ (cf. \cite{BLLM25}).  
It is thus natural to ask whether, for any polynomial $Q$ which has a simple
critical point $c_Q$, there exists an involution J which fixes $c_Q$ and yields a
correspondence $J\circ \mathrm{Cov}_0^Q$, which is mating between some parabolic
rational map and the Hecke group.
\begin{figure}[h!] 
\captionsetup{width=0.96\linewidth}
\includegraphics[width=0.96\linewidth]{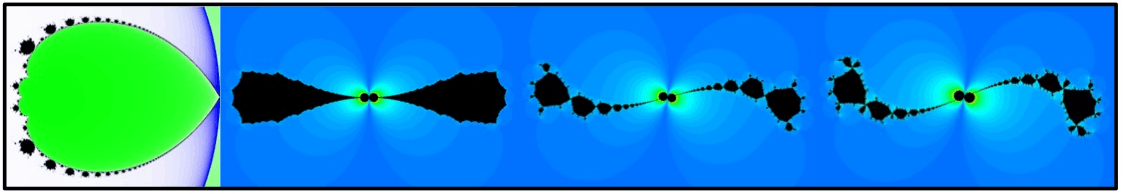}
\includegraphics[width=0.96\linewidth]{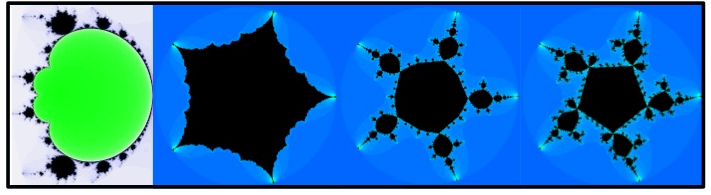}
    \caption{Top: Modular Multibrot set  and some limit sets. Bottom: Parabolic Multibrot set and corresponding Julia sets.}
  \label{mult}
\end{figure}
In \cite{BLLS}, a first step towards addressing this deep problem is taken in the holomorphic setting. There it is shown that, for every degree $d+1$ polynomial $Q$ with two critical points in $\C$ (one of which, say $c_Q$, is simple), if $J_b$ is an involution fixing $c_Q$ and $b \in \C$, then the family of correspondences $\cF_b := J_b \circ Cov^Q_0$
contains matings between the Hecke group $\mathcal{H}_{d+1}$ and all bicritical parabolic rational maps $R$ satisfying the hypothesis of Theorem~\ref{par_rat_hecke_thm} (also in the case when $R$ has a disconnected Julia set).  
Moreover, the authors establish the existence of a dynamically natural bijection between the connectedness locus for the family $\cF_b$, termed the \emph{Modular Multibrot}, and the connectedness locus for the family of bicritical parabolic rational maps $R$, referred to as the \emph{Parabolic Multibrot}. They conjecture that this bijection is in fact a homeomorphism (see Figure~\ref{mult} for numerical evidence). 

The above setting can be generalised to the case where $Q$ is a Belyi polynomial \cite{BLLS3}. We refer the reader to \cite{LMM25d} where the general Belyi case was settled in the antiholomorphic world.

\section{Where Teichm{\"u}ller spaces meet connectedness loci.}\label{para_unif_subsec}
The results discussed in \S\ref{gen_matings_sec} shed light on the co-existence of parameter spaces of polynomials and Kleinian groups in spaces of correspondences.
\smallskip

\noindent\textbf{Product structure in spaces of correspondences.}
In \cite[Appendix~A]{LLM24}, `mating continuity' arguments were used to show  that for $\Sigma\in\mathfrak{F}$ (see \S\ref{holo_gen_mating_subsec}), the space $\mathrm{Teich}(\Sigma)\times\mathcal{C}_{d_\Sigma,fr}$ (where $\mathcal{C}_{k,fr}$ is the space of degree $k$ periodically repelling, finitely renormalizable polynomials with connected Julia set) naturally embeds into a space of correspondences. The image of this embedding also lives in the space of degenerate p-l maps, where the product structure comes from varying the internal and the external map independently. Another important topological embedding is that of $\mathrm{Teich}(\Sigma)\times \cH_{d_\Sigma}$ (where $\cH_k$ is the main hyperbolic component of degree $k$ polynomials) into a space of bi-degree $d_\Sigma:d_\Sigma$ algebraic correspondences parametrized by degree $(d_\Sigma+1)$ rational maps; this is proved using quasiconformal surgery methods \cite[\S 6]{LMM25a}. The image of this embedding can be regarded as a \emph{hybrid simultaneous uniformization locus} lying halfway between quasi-Fuchsian and quasi-Blaschke spaces (cf. \cite{Ber60,BF14}). The image of a slice $\mathrm{Teich}(\Sigma)\times \{P\}$ under this embedding, where $P\in\cH_{d_\Sigma}$, is called the \emph{Bers slice} of $\Sigma$ passing through $P$, and is denoted by $\mathfrak{B}(P)$.
The slice $\mathfrak{B}(z^{d_\Sigma})$ is of particular interest: it embeds holomorphically in the space of correspondences and admits an explicit description \cite[\S 7]{MM25}.

In \cite{BH00}, the moduli space of a $2$ parameter family of correspondences was shown to be in bijection with $\left(\Int{\mathcal D}\right)/\iota \times \cM$, where $\iota$ is an involution (and $\Int{\mathcal D}$ is the space of discrete faithful representations of $C_2*C_3$ in $\PSL(2,\C)$ with Cantor limit set).
In \cite{BH07} matings between $z^2$ and every `cusp' group on $\partial \mathcal{D}$
are realized by pinching techniques.
The results in \cite{BH00} are refined and generalized to degree $d$ polynomials and $C_2*C_{d+1}$ in \cite{Lau24}, where the mating operation is shown to be analytic on the interior and continuous at quasiconformally rigid parameters. Thus, Ratis Laude constructed topological products of quasi-Fuchsian deformations of certain groups and subsets of polynomial connectedness loci in spaces of correspondences \cite{Lau24}.
\smallskip

\noindent\textbf{Boundedness results and Bers boundary correspondences.}
Using rescaling limit techniques, it was proved in \cite{LMM25a} that the Bers slices $\mathfrak{B}(P)$ are pre-compact in the space of bi-degree $d_\Sigma:d_\Sigma$ correspondences on $\widehat{\C}$, and the compactifications are naturally homeomorphic when $\mathrm{dim}_\C\mathrm{Teich(\Sigma)=1}$. This result, which is an analog of pre-compactness of classical Bers slices \cite{Ber70}, opens the door for a detailed study of \emph{Bers boundary correspondences}. 
\smallskip

\noindent\textbf{Teichm{\"u}ller, Blaschke, and Hurwitz.}
The ambient space of the correspondences constructed in \cite{MV25} is an appropriate \emph{Hurwitz space} of pairs $(\Sigma,\cR)$, where $\Sigma$ runs over the moduli space $\mathcal{M}_g$ of genus $g$ surfaces, and $\cR$ ranges over meromorphic maps defined on such surfaces. The mating construction yields a dynamically natural injection of products of finitely many Teichm{\"u}ller and Blaschke spaces into the said Hurwitz space (see \cite{MV25}).

\section{When Julia sets look like limit sets.}\label{limit_julia_homeo_sec}
In this section, we summarize recent works on the existence of equivariant homeomorphisms between various Julia and limit sets, and their parameter space implications.
\subsection{A dynamical link between Julia and limit sets.}\label{dyn_dict_subsec}
\begin{figure}[htbp]
\captionsetup{width=0.96\linewidth}
  \centering
  \includegraphics[width=0.99\linewidth]{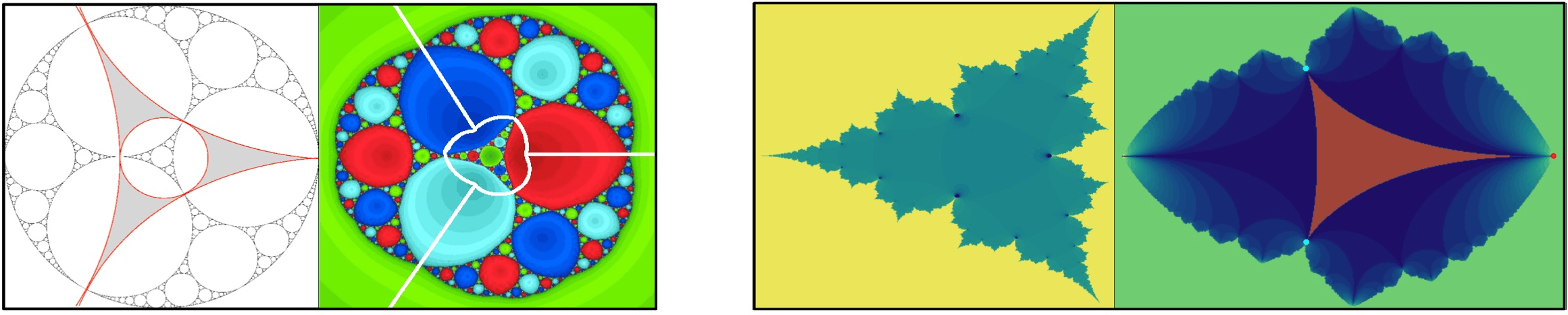}
  \caption{Left: The topological realization of the Apollonian gasket as a Julia set. Right: Two conformally removable non-quasicircle Jordan curves, one with two-sided cusps and the other with one-sided cusps, arising as limit sets of matings.}
  \label{apollo_fig}
\end{figure}

The Nielsen map $\pmb{\cN}_d$ of the ideal $(d+1)$-gon reflection group $\pmb{G}_d$ (see \S\ref{sigma_d_subsec}) that act as a bridge between anti-polynomials and reflection groups in the mating theory, also gives rise to precise connections between suitable Julia and limit sets. This was first exploited in \cite{LLMM23b} to prove that the Apollonian gasket, which is the limit set of a reflection group, is homeomorphic to the Julia set of a critically fixed cubic anti-rational map such that the homeomorphism conjugates the action of the Nielsen map to that of the anti-rational map (see Figure~\ref{apollo_fig}). In \cite{LMM22}, limit sets of necklace groups (see \S\ref{sigma_d_subsec}) were shown to be naturally homeomorphic to Julia sets of critically fixed anti-polynomials. The key idea in all these constructions is to organize the circular reflections generating a reflection group into a single map; i.e., the Nielsen map. The action of this Nielsen map, restricted to an `invariant component' of the ordinary set, is topologically conjugate to $\pmb{\cN}_k$, for some $k\geq 2$. One introduces critical points by replacing this map $\pmb{\cN}_k$ with the power map $\overline{z}^k$, thus producing a critically fixed branched cover of $\mathbb{S}^2$. Finally, an application of W. Thurston's realization theorem (see \cite{DH93,PL98}) upgrades the branched cover to a critically fixed anti-rational map. Conversely, one can start with a critically fixed anti-rational map, and replace the $\overline{z}^k$-dynamics on an invariant Fatou component with that of the Nielsen map $\pmb{\cN}_k$ via David surgery, producing the desired reflection group \cite[\S 8]{LMMN25} (see \cite[\S 8.1]{LM23} for detailed discussions). The most general version of this surgery procedure led to the following result. By definition, a \emph{kissing reflection group} is the group generated by reflections in the circles of a finite (connected) circle packing.

\begin{theorem}\cite{LLM22a,LMMN25}
There exists a natural bijective correspondence between:
\begin{enumerate}
\item $\{2$-connected, simple, plane graphs $\mathscr{G}$ with $d+1$ vertices up to isomorphism of plane graphs$\}$,

\item $\{$Kissing reflection groups $G=G_{\mathscr{G}}$ of rank $d+1$ with connected limit set up to QC conjugacy$\}$, and

\item  $\{$Critically fixed anti-rational maps $R=R_{\mathscr{G}}$ of degree $d$ up to M{\"o}bius conjugacy$\}$.
\end{enumerate}

\noindent Here, $G_\mathscr{G}$ is a group generated by reflections in the circles of a circle packing having $\mathscr{G}$ as its contact graph. The Tischler graph of $R_\mathscr{G}$ (which is the union of all invariant internal rays in the critical Fatou components) is the planar dual of $\mathscr{G}$. Moreover, there exists a David homeomorphism of $\widehat{\C}$ that conjugates $R_\mathscr{G}\vert_{\mathcal{J}(R_\mathscr{G})}$ to $\cN_{G_\mathscr{G}}\vert_{\Lambda(G_\mathscr{G})}$.
\end{theorem}

This bijection respects several dynamical aspects. In particular, the sub-class of necklace groups correspond to critically fixed anti-polynomial (both have completely invariant, simply connected domains in their normality set), and groups on quasi-Fuchsian boundaries (i.e., matings of necklace groups) correspond to anti-rational maps that are matings of two anti-polynomials. These properties can also be detected from the graph $\mathscr{G}$.

A characterization of critically fixed anti-rational maps was obtained independently by Geyer \cite{Gey20}. In the holomorphic setting, Julia sets of polynomials that are naturally homeomorphic to limit sets of certain Kleinian punctured sphere groups were constructed in \cite{MM23a}, using Bowen-Series and higher Bowen-Series maps as replacements of Nielsen maps. Recently, it was shown in \cite{LMM25c} that limit sets of geometrically finite Kleinian groups lying on Bers boundaries of closed surfaces are quasiconformally homeomorphic to the Julia set of the \emph{fat Basilica} polynomial $z^2-3/4$, using so-called \emph{Basilica Bowen-Series maps}.

\subsection{Deformation space analogies.}
The link between kissing groups and  anti-rational maps (\S\ref{dyn_dict_subsec}) suggests parallels between their parameter spaces. The quasiconformal deformation space of a kissing group serves as its parameter space, where the parabolic elements stay parabolic. These parabolics correspond to the repelling fixed points of the corresponding anti-rational map. To reflect this `persistently parabolic' feature, one restricts to a \emph{pared deformation space} inside the hyperbolic component of a critically fixed anti-rational map, where repelling fixed points have uniformly bounded multipliers in suitable uniformizing models.
With these notions of parameter spaces, the following analogies were unearthed in \cite{LLM22b} (see \cite[\S 8.2]{LM23} for a detailed survey).
\smallskip

\noindent\textbf{Boundedness results.} Basic results from $2$ and $3$-dimensional topology show that $\Lambda(G_\mathscr{G})$ and $\mathcal{J}(R_\mathscr{G})$ are \emph{gaskets} iff the $3$-manifold $\mathbb{H}^3/\widetilde{G}_\mathscr{G}$ (where $\widetilde{G}_\mathscr{G}$ is the index $2$ Kleinian subgroup of $G_\mathscr{G}$) is acylindrical, which is further equivalent to $\mathscr{G}$ being polyhedral or $3$-connected. By \cite{Thu86a}, this is precisely the condition for boundedness of the deformation space of $G_\mathscr{G}$. Using rescaling limits and degenerations to limiting trees, it was shown in \cite{LLM22b} that $3$-connectedness of $\mathscr{G}$ is also equivalent to boundedness of the pared deformation space of $R_\mathscr{G}$. 
\smallskip

\noindent\textbf{Bifurcations and global topology} The closure of the deformation space of $G_\mathscr{G}$ is stratified into deformation spaces of kissing groups whose underlying circle packings have more tangencies. This can be captured by a graph-theoretic relation called `domination'. By realizing appropriate degenerations, it was proved in \cite{LLM22b} that the same graph-theoretic condition also determines `bifurcations' of critically fixed hyperbolic components from the hyperbolic component of $R_\mathscr{G}$. 
This above bifurcation structure implies, among other things, that the union of the closures of pared
deformation spaces of critically fixed anti-rational maps of a given degree is connected. We refer the reader to \cite[Theorem~1.4]{LLM22b} for a result about the topological complexity of this union that mirrors a similar result on the moduli space of circle packings proved in \cite{HT}.

\section{Applications near and far.}\label{applications_sec} The methods developed to explore links between rational dynamics and Kleinian groups have applications in other areas of mathematics. In this section, we present a selection of them.
\smallskip

\noindent\textbf{Conformal removability.} A compact set $E\subset \widehat{\C}$ is \textit{conformally removable} if every homeomorphism $f\colon \widehat{\C}\to \widehat{\C}$ that is conformal on $\widehat{\C} \setminus E$ is a M{\"o}bius map (cf. \cite{You15}). Conformal removability of limit and Julia sets plays an important role in dynamics and geometric function theory. Removability of connected Julia sets of semi-hyperbolic polynomials follows from the John property of their basins of infinity \cite{CJY94,Jon95}. The presence of parabolics/cusps makes the situation more subtle. Using David surgery techniques, it was proved in \cite[\S 9]{LMMN25} that connected Julia sets of geometrically finite polynomials and limit sets of necklace reflection groups (which are `cuspidal fractals')  are conformally removable. The quasiconformal compatibility results of \cite{LMM25c} further imply that limit sets of geometrically finite Kleinian groups lying on Bers boundaries of closed surfaces are conformally removable.
\smallskip

\noindent\textbf{Conformal welding.} A homeomorphism $h\colon \mathbb S^1 \to \mathbb S^1$ is called a \emph{welding homeomorphism} if the Riemann maps of the complementary components of a Jordan curve $J$ differ by $h$. Such a Jordan curve $J$ is called a \emph{welding curve} for $h$. It is a classical fact that quasisymmetric circle homeomorphisms are welding maps for quasicircles, see \cite{Ham91,Bis07} for more examples of welding homeomorphisms. Matings of polynomials with groups lead to certain `equivariant welding problems', which gave rise to a general realization theorem for conjugacies between suitable piecewise analytic, expansive circle coverings as welding homeomorphisms (see Figure~\ref{apollo_fig}, \cite[\S 5]{LMMN25}).
\smallskip

\noindent\textbf{Quasisymmetry groups.} Quasiconformally equivalent fractals have isomorphic self-quasisymmetry groups, makings the quasisymmetry group of a fractal an important invariant in metric geometry. Quasisymmetry groups of some carpet Julia/limit sets and polynomial Julia sets were studied previously in \cite{BLM16,LM18} (also see \cite{All23,BF25}). The results arising from the theme of this survey brought to light certain limitations of the use of quasisymmetry groups as quasiconformal invariants; indeed, according to \cite{LLMM23b} certain homeomorphic (gasket) Julia and limit sets have isomorphic quasisymmetry groups, but they were shown to be quasiconformally nonequivalent in \cite{LZ25}.  
\smallskip

\noindent\textbf{Topology and singularities of quadrature domains and algebraic droplets.} In \cite{LM16}, bounds on connectivity of quadrature domains were given using quasiconformal surgery of Schwarz reflections and Hele-Shaw flow techniques, and these bounds were shown to be sharp in \cite{LM14}. In \cite{MR25}, alternative and non-perturbative arguments involving the dynamics of Schwarz reflections were used to combine the bounds of \cite{LM16} with linear bounds on the number of weighted singular points of a quadrature domain. This, in particular, improved the previously known quadratic upper bound on the singularities of a quadratic domains \cite{Gus83} to a linear one.
\smallskip

\noindent\textbf{Univalent rational maps and extremal problems.} Univalent restrictions of polynomials/rational maps is a classical area in complex analysis \cite[\S 7.4]{SS02}. The dynamics of Schwarz reflections and their quasiconformal deformations can be used profitably to study such problems of analytic origin. We refer the reader to \cite{LMM21,LMM25d} (cf. \cite[\S 12.4]{LM23}) for recent progress in this area obtained by employing dynamical tools.  
Rational maps in classes $S, \Sigma$ of univalent maps and their extremal points have also received considerable attention in the study of coefficient optimization \cite{Suf72,Dur83,SS02} and quadrature/liquid domains \cite{ADPZ20,LM14}. A dynamical theory of Schwarz reflections arising from $\Sigma_d^*\subset\Sigma$ and $S_d^*\subset S$ (see \S\ref{sigma_d_subsec}) led to a classification of extremal points of these spaces of univalent functions. Connections between these extremal points and maximal cusp necklace groups were unearthed \cite{LMM21,LMMN25}.
\smallskip

\noindent\textbf{Extension of Bers Simultaneous Uniformization Theorem.} The classical Bers Simultaneous Uniformization Theorem allows one to combine a pair of homeomorphic surfaces, carrying different conformal structures, into a single uniformizing Kleinian group \cite{Ber60}. The introduction of factor Bowen-Series maps renders more flexibility to this construction. In \cite{MM24}, ideas from the mating program of \S\ref{gen_matings_sec} were adapted for a purely surface-theoretic setting, and algebraic correspondences on $\widehat{\C}$ were constructed as combinations of certain pairs of Fuchsian groups uniformizing \emph{non-homeomorphic} genus zero orbifolds. As an application, new complex-analytic realizations of Teichm{\"u}ller spaces of genus zero orbifolds were found in parameter spaces of algebraic correspondences. 
\smallskip

\noindent\textbf{Orbit equivalence.} In \cite[\S 8]{MM23a}, Bowen-Series maps of Fuchsian punctured sphere groups were used to demonstrate failure of topological orbit equivalence rigidity for Fuchsian actions on $\mathbb{S}^1$. A partial classification of piecewise M{\"o}bius circle coverings that are orbit equivalent to arbitrary Fuchsian groups was given in terms of (higher) Bowen-Series maps in \cite[\S 5, \S 6]{MM23a},
On the other hand, by introducing additional break-points and accelerating the dynamics of Bowen-Series maps, a new class of piecewise M{\"o}bius circle coverings, generically orbit equivalent to arbitrary (non-rigid) Fuchsian groups, have been constructed in \cite{LMM25c}.

\section{Related works.}\label{connections_sec}
Holomorphic correspondences were studied from an ergodic-theoretic point of view by various authors; e.g., constructions of invariant, ergodic measures, statements about their supports, and equidistribution results can be found for 
\begin{enumerate}[leftmargin=8mm]
    \item modular correspondences in \cite{CU03,Din13},
    \item correspondences of unequal bi-degree in \cite{DS06,BS16},
    \item \emph{weakly modular} correspondences of equal bi-degree in \cite{DKW20},
    \item the correspondences $\cF_a$ treated in \S\ref{BP_sec} and their generalizations in \cite{Par24,Par23}, and
    \item correspondences arising from holomorphic semigroups in \cite{BS17,Lon22a,Lon22b}.
\end{enumerate}

For analytic applications of equidistribution results for correspondences, see \cite{Hem24}. A theory of thermodynamic formalism for correspondences involving pressure and variational principle was developed in \cite{LLZ23,Wan25}, and conformal measures for several correspondences (including many matings) were constructed in \cite{HLL24}. For results about Hausdorff dimension, holomorphic motions of Julia sets, and phase--parameter similarity near Misiurewicz parameters for correspondences of the form $z^r + c$, $r\in\mathbb{Q}$, see \cite{SS17,Siq23,Siq25a,Siq25b}. 
Correspondences were investigated through the lens of arithmetic dynamics in \cite{Ing17,Ing19}.
Various correspondences arising from Thurston theory were studied in \cite{Ram18,BDP24}. 
Deformations of rigid Kleinian groups in the wider category of correspondences were constructed in \cite{BP99}. A Bers slice of matings between the modular group and its correspondence deformations is constructed and analysed in \cite{B10b}.
An example of infinite bi-degree correspondence realizing the mating of a transcendental map with a Fuchsian group was given in \cite{BF08}.

Quadrature domains also appear in the study of zeroes of Gaussian Entire Functions \cite{NW24}. The deltoid reflection map of \S\ref{deltoid_subsec} emerged in certain limiting study of dessin d'enfants in \cite{ILRS23}. 
For results on renormalization and boundedness in parameter spaces of circle packings, and uniformizations of gasket Julia sets, see \cite{LZ23,LZ24,LN24b}.
The David extension theorem of \cite{LMMN25} is based on certain general sufficient conditions obtained in \cite{CCH96,Zak08}. A sharper sufficient condition for David extensions of circle homeomorphisms was later proved in~\cite{KN22}. For quasiconformal equivalence of dynamically arising infinitely cusped Jordan curves (including limit sets of various conformal matings)  and the welding curve of the Minkowski $?$-function, see \cite{McM25}. We refer the reader to \cite{Los14,AKU22a,AKU22b,AJLM23,MM23b,AKU25} for various results on entropy of Bowen-Series type maps, and connections with Gromov boundaries of hyperbolic groups and reduction theory for Fuchsian groups.

\section{Further directions.}\label{new_directions_sec}
\begin{enumerate}[leftmargin=8mm]
\item Classify  correspondence realizing matings
of rational maps and Kleinian groups; e.g., of the form $J\circ \mathrm{Cov}_0^R$ (cf. \S\ref{inverse_subsec}).
\item Develop a dynamical theory for the family of correspondences $J\circ \mathrm{Cov}_0^R$ (cf. \cite{BL22}), including their behavior outside the \emph{mating locus} and the \emph{discreteness locus}.
\item Construct matings of arbitrary Fuchsian groups with polynomials, following the program of \S\ref{gen_matings_sec}.
\item Extend the combination program to mate transcendental maps with (possibly infinitely generated) Fuchsian groups (cf. \cite{BF08}). 
\item Study product structures, degeneration phenomena, double limit theorems in parameter spaces of correspondences and their connections with Hurwitz spaces (see \cite[\S 2.6]{LM23}, \cite[\S 1.4]{MV25} for precise questions).
\item Mate polynomials $P$ (connected $K(P)$) with Bers boundary groups as correspondences, possibly with Peano curve limit sets.
\item Investigate the dynamics of Bers boundary correspondences (\S\ref{para_unif_subsec}) and the topological structure of such Bers boundaries (see \cite[\S 1.3]{LMM25a}, \cite[\S 4.5]{MM24} for precise questions).
\item Study continuity of the mating operation going to the boundary of Teichm{\"u}ller spaces (cf. \cite{BH07}). 
\item Develop a bifurcation theory for the Mandelbrot sets of correspondences studied in \cite{Ing17,Siq25a}.
\item Study spaces of degenerate p-l maps as limits of p-l maps. 
\item Combine different correspondences within the space of correspondences (cf. \cite{B10}).
\item Design an analog of Lyubich-Minsky laminations for correspondences \cite{LM97}, and use it to study equidistribution and rigidity in parameter spaces.
\item Construct Weil–Petersson/pressure metrics for the hybrid simultaneous uniformization loci of \S\ref{para_unif_subsec}. 
\item Classify fractals (e.g., Jordan curves) arising as limit sets of matings of hyperbolic/parabolic maps 
up to quasiconformal equivalence (cf. \cite{McM25}). 
\item Study Galois equidistribution in parameter spaces of correspondences. 
\end{enumerate}

\section*{Acknowledgments.}
We are grateful to Shaun Bullett, Mikhail Lyubich, Nikolai Makarov, Yusheng Luo, and Mahan Mj for sharing their ideas and insights with us. S. M. also thanks Kirill Lazebnik, Seung-Yeop Lee, Russell Lodge, Jacob Mazor, Sergiy Merenkov, and Dimitrios Ntalampekos for many fruitful collaborations.

\end{document}